\documentclass{amsart}
\usepackage{amssymb,amsmath, amsthm,latexsym}
\usepackage{graphics}
\usepackage{amscd}
\usepackage{graphics}
\newcommand{\cal}[1]{\mathcal{#1}}
\theoremstyle{plain}
\newtheorem*{theo}{Theorem}

\newtheorem{lemma}{Lemma}[section]
\newtheorem{theorem}[lemma]{Theorem}
\newtheorem{proposition}[lemma]{Proposition} 
\newtheorem{corollary}[lemma]{Corollary}
\theoremstyle{definition}
\newtheorem{definition}[lemma]{Definition}
\newtheorem{remark}[lemma]{Remark}
\newtheorem{example}[lemma]{Example}
\parskip=\bigskipamount

\let\egthree=\phi
\let\phi=\varphi
\let\varphi=\egthree

\begin{document}
\title{Lines of Minima in Outer space}
\author{Ursula Hamenst\"adt}
\date{June 10, 2013}
\begin{abstract} 
We define lines of minima in the thick part of Outer space
for the free group $F_n$ with $n\geq 3$ generators.
We show that these lines of minima are contracting
for the Lipschitz metric. Every fully irreducible 
outer automorphism of $F_n$ defines such a line of minima.
Now let $\Gamma$ be a subgroup
of the outer automorphism group of $F_n$ which is not
virtually abelian. We obtain that 
if $\Gamma$ contains at least one
fully irreducible element then for every 
$p\in (1,\infty)$ the
second bounded cohomology group 
$H_b^2(\Gamma,\ell^p(\Gamma))$ is infinite dimensional.
\end{abstract}

\maketitle

\section{Introduction}

There are many resemblances between the \emph{extended mapping
class group} of a closed oriented surface $S$, i.e the outer automorphism
group of the fundamental group of $S$,  
and the outer automorphism group
${\rm Out}(F_n)$ of a free group $F_n$ with $n\geq 2$
generators. Most notably, ${\rm Out}(F_2)$ is just the
extended mapping class group $GL(2,\mathbb{Z})$ 
of a torus. However, while in recent years the attempt to 
understand the mapping class group via the geometry of
spaces on which it acts lead to a considerable gain of knowledge
of the mapping class group, 
so far this approach has not been
carried out successfully for ${\rm Out}(F_n)$.

The group ${\rm Out}(F_n)$ acts properly
discontinuously on \emph{Outer space} ${\rm CV}(F_n)$.
This space consists of projective classes of 
marked metric graphs
with fundamental group $F_n$ and 
can be viewed as an equivalent of Teichm\"uller space
for a closed surface $S$ of higher genus.
Teichm\"uller space admits several natural and
quite well understood metrics which are invariant under the
action of the extended mapping class group.
The best known such metrics are the 
\emph{Teichm\"uller metric} and the \emph{Weil-Petersson metric}. 
In contrast, up to date there is no
good geometric theory 
of Outer space. Only very recently 
Francaviglia and Martino \cite{FM08} studied in a systematic
way a natural ${\rm Out}(F_n)$-invariant
metric on Outer space. 
This metric is the symmetrization of a non-symmetric
geodesic metric, the so-called \emph{Lipschitz metric}.
The symmetric metric is not
geodesic, and its analogue for Teichm\"uller space, the
\emph{Thurston metric}, also turned out to be 
harder to understand than the Teichm\"uller metric
and the Weil-Petersson metric.

The viewpoint we take in this work is motivated by
a slightly different approach to
the geometry of Teichm\"uller space.
Namely, \emph{lines of minima} in Teichm\"uller space
for a closed surface $S$ of higher genus were defined and
investigated by Kerckhoff \cite{Ke92}. 
As for geodesics for the Teichm\"uller metric, 
such a line of minima
is determined by two projective measured geodesic laminations 
which jointly fill up $S$. 
A line of minima uniformly fellow-travels its corresponding 
Teichm\"uller geodesic provided that this Teichm\"uller geodesic
entirely remains in the \emph{thick} part of Teichm\"uller space.
There is also a very good control in the thin part of 
Teichm\"uller space though the uniform fellow traveller
property is violated \cite{CRS08}. 

More precisely,
for a number $\epsilon >0$, the $\epsilon$-thick part 
${\cal T}(S)_\epsilon$ of 
Teichm\"uller space for a closed oriented 
surface $S$ of higher genus is the set of all hyperbolic 
metrics on $S$ whose systole (i.e. the length of a shortest closed
geodesic) is at least $\epsilon$.
Teichm\"uller geodesics $\gamma$ entirely
contained in ${\cal T}(S)_\epsilon$, 
and hence their corresponding lines of minima, have
a \emph{uniform contraction property}:  
If $B$ is a closed metric ball in Teichm\"uller space 
for the Teichm\"uller metric or the Weil-Petersson metric
which is disjoint from $\gamma$ 
then the diameter of 
a shortest distance projection of $B$ into $\gamma$
is bounded from above by a constant only depending on
$\epsilon$ \cite{Mi96,BFu09}.
 
Our main goal is to define such lines of minima
for Outer space. We show that these lines of minima are
uniform coarse geodesics for the
symmetrized Lipschitz metric, and these
coarse geodesics have the uniform contraction
property (Corollary \ref{contraction2}).
We also observe that for 
every fully irreducible element 
$\phi\in {\rm Out}(F_n)$ there is such 
a line of minima which is $\phi$-invariant.
This recovers a recent result of 
Algom-Kfir \cite{AK08}
who showed that axes for fully irreducible elements as defined by
Handel and Mosher \cite{HM06} have the contraction property.

As an application, we use the tools
developped in
\cite{H08} to show

\begin{theo}\label{mainapplication}
Let $\Gamma<{\rm Out}(F_n)$ be a subgroup which is
not virtually abelian and which contains at least
one fully irreducible element. Then for every $p>1$ the
second bounded cohomology group
$H_b^2(\Gamma,\ell^p(\Gamma))$ is infinite dimensional.
\end{theo}

Earlier Bestvina and Feighn \cite{BF10} showed that
a subgroup $\Gamma$ as in the theorem has 
nontrivial second bounded cohomology with real coefficients. 
This is also immediate from our approach.

All constructions in this paper are equally valid for the
action of the mapping class group on Teichm\"uller space. 
This leads for example to a new proof of the main result of 
\cite{Mi96} avoiding completely the explicit use of Teichm\"uller theory.
However, in this case our more abstract
approach does not have any obvious
advantage over the original arguments. The analogue of 
our main theorem for mapping class groups was derived with
a different method in \cite{H08b}.

\section{Measured laminations and trees}

In this section we introduce currents, trees and
measured laminations for the free group
$F_n$ of rank $n\geq 3$. 
We single out an ${\rm Out}(F_n)$-invariant subset
of the space of measured laminations which
is used for the construction of lines of minima
in the later sections. We continue to use the
notations from the introduction.

The Cayley graph of $F_n$ with respect to
a fixed standard symmetric generating set is a
regular simplicial tree which can be compactified by adding
the \emph{Gromov boundary} $\partial F_n$. This boundary
is a compact totally disconnected topological space
on which $F_n$ acts as a group of 
homeomorphisms. It does not depend on the generating set
up to $F_n$-equivariant homeomorphism.
Every element $w\not=e\in F_n$ acts on
$\partial F_n$ with \emph{north-south dynamics}.
This means that $w$ fixes precisely two points 
$a_+,a_-\in \partial F_n$, 
and for every neighborhood $U_+$ of $a_+$,
 $U_-$ of $a_-$ there is some $k>0$ such that
$w^k(\partial F_n-U_-)\subset U_+$ and
$w^{-k}(\partial F_n-U_+)\subset U_-$. 

A \emph{geodesic current} for $F_n$ is a locally finite
Borel measure on 
\[\partial F_n\times \partial F_n-\Delta=\partial ^2(F_n)\]
(where $\Delta$ denotes the diagonal in $\partial F_n\times
\partial F_n$)
which is invariant under the action of $F_n$ and under
the flip $\iota: \partial^2( F_n)\to \partial^2(F_n)$ 
exchanging the two factors. 
The space ${\rm Curr}(F_n)$ of all geodesic 
currents equipped with the
weak$^*$-topology is a locally compact 
topological space which can
be projectivized to the compact space ${\cal P}{\rm Curr}(F_n)$
of \emph{projective
currents}. The \emph{outer automorphism group} ${\rm Out}(F_n)$
of $F_n$ naturally acts on ${\rm Curr}(F_n)$ and on
${\cal P}{\rm Curr}(F_n)$ as a group of homeomorphisms.

If $w\not=e\in F_n$ is any indivisible element, 
i.e. an element which is not of the form
$w=v^k$ for some $k\geq 2$, then 
the set of all pairs of
fixed points in $\partial F_n$ 
of all elements of $F_n$ which are conjugate to $w$
is a discrete $F_n$-invariant flip invariant   
subset of $\partial F_n\times \partial F_n-\Delta$. Thus 
the sum of the Dirac measures supported at these 
pairs of fixed points
defines a geodesic current
which we call \emph{induced} by the conjugacy class $[w]$ of $w$. 
If $w=v^k$ for some $k\geq 2$ and some 
indivisible element $v\in F_n$ then we define
the geodesic current induced by the conjugacy class $[w]$ to 
be the $k$-fold multiple of the geodesic current
induced by  
the conjugacy class $[v]$ of $v$. Define a
\emph{weighted induced current} to be a geodesic 
current which is obtained
by multiplying a geodesic current induced by a conjugacy
class in $[F_n]$ by a positive weight. The set of weighted induced
currents is invariant
under the action of ${\rm Out}(F_n)$.

An element $w\not=e\in F_n$ is \emph{primitive} if
it belongs to some basis, i.e. if there is a decomposition
of $F_n$ into a free product of the form
$<w>*H$ where $<w>$ is the infinite cyclic subgroup of
$F_n$ generated by $w$ and where
$H$ is a free subgroup of $F_n$.
A conjugacy class in $F_n$ is \emph{primitive} if one 
(and hence each) of its elements is primitive.
The set of primitive elements is invariant under
the action of the full automorphism group ${\rm Aut}(F_n)$ of
$F_n$ and hence
${\rm Out}(F_n)$ naturally acts on the set of all
primitive conjugacy classes.

\begin{definition}\label{measlam}
The space ${\cal M\cal L}(F_n)$
of \emph{measured laminations} is the closure
in ${\rm Curr}(F_n)$ of the 
set of all currents which are weighted
induced currents of primitive conjugacy classes. 
\end{definition}

The projectivization 
${\cal P\cal M\cal L}(F_n)$
of ${\cal M\cal L}(F_n)$, equipped with the
weak$^*$-topology, is compact and invariant under
the action of ${\rm Out}(F_n)$.
Theorem B of
\cite{KL07a} shows that 
${\cal P\cal M\cal L}(F_n)$ 
is the unique smallest non-empty
closed ${\rm Out}(F_n)$-invariant subset of 
${\cal P}{\rm Curr}(F_n)$.

Martin (Theorem 17 of \cite{Ma95}) characterizes 
projective measured laminations
which are induced by a conjugacy class 
in $F_n$ as follows.

\begin{proposition}\label{properfreefactor} 
A projective current induced by a
conjugacy class $[\alpha]$ in $F_n$ is contained in 
${\cal P\cal M\cal L}(F_n)$ if and only if
each element of $[\alpha]$ is contained in a proper
free factor of $F_n$. 
\end{proposition}

A closed non-empty $F_n$-invariant
subset of $\partial^2(F_n)$ which is moreover invariant
under the flip $\iota$ is called a
\emph{topological lamination}.
The space ${\cal L}$ 
of all topological laminations 
can be equipped with the \emph{Chabauty topology}.  
With respect to this topology, ${\cal L}$ is compact.
The group ${\rm Out}(F_n)$ acts on ${\cal L}$ 
as a group of homeomorphisms.
Every nontrivial element $w\not=e\in F_n$ defines a 
point 
$ [[w]]\in {\cal L}$ which is just the 
set of all pairs of fixed points of 
all elements of $F_n$ which are conjugate to $w$. In other
words, $[[w]]$ is the support of the current induced by
the conjugacy class of $w$.
Topological laminations of this form are 
called \emph{rational}.
The support of a measured lamination is a topological
lamination.

We call a (projective)
geodesic current supported in 
a topological lamination $L$ a \emph{(projective) transverse measure}
for $L$. If $L$ is the topological
lamination defined by the conjugacy class
of a primitive element $w\in F_n$ then $L$ admits a
unique projective transverse measure. This measure is just
the projective measured lamination induced by the conjugacy class
of $w$.
Note that if $L_i\to L$ in ${\cal L}$ and
if $\zeta_i$ is a projective geodesic
current supported in $L_i$ then up to
passing to a subsequence, the projective geodesic
currents
$\zeta_i$ converge in ${\cal P}{\rm Curr}(F_n)$ 
to a projective geodesic current
supported in $L$ (we refer to \cite{CHL07b} for a more
precise discussion).

\begin{remark} The  
definition of a topological lamination
does not correspond to the definition of a geodesic lamination
for closed surfaces. The correct analogue of a lamination
in the surface case is a closed $F_n$-invariant subset
of $\partial^2(F_n)$ which is contained in the Chabauty closure
of those closed $F_n$-invariant subsets of $\partial^2(F_n)$ 
which consist of pairs of fixed points of elements in some primitive
conjugacy class. 
\end{remark}

Let $cv(F_n)$ be the space of all minimal free and
discrete isometric actions of $F_n$ on $\mathbb{R}$-trees.
Two such actions of $F_n$ on $\mathbb{R}$-trees 
$T$ and $T^\prime$ are identified in
$cv(F_n)$ if there exists an $F_n$-equivariant isometry
between $T$ and $T^\prime$. 
The quotient of a tree $T\in cv(F_n)$ 
under the action of $F_n$ is a
finite metric graph $T/F_n$ without vertices of valence one or two 
whose fundamental group is marked isomorphic to $F_n$.
The space $cv(F_n)$ admits a natural locally compact
metrizable ${\rm Out}(F_n)$-invariant topology.

The \emph{boundary} $\partial cv(F_n)$ of $cv(F_n)$ consists
of all minimal \emph{very small} isometric actions 
of $F_n$ on $\mathbb{R}$-trees
which either are non-simplicial or which are not free.
Here an action is very small if and only if
every nontrivial arc stabilizer is maximal cyclic and if
tripod stabilizers are trivial. Again, any two such
actions define the same point in $\partial cv(F_n)$ if
there exists an $F_n$-equivariant isometry between them.
The union
$\overline{cv(F_n)}=cv(F_n)\cup \partial cv(F_n)$ has a natural
${\rm Out}(F_n)$-invariant topology such that 
$cv(F_n)\subset \overline{cv(F_n)}$ is open and dense.

\emph{Outer space} ${\rm CV}(F_n)$ is 
the projectivization of 
$cv(F_n)$. If we define two trees $T,T^\prime$ to be
equivalent if there is an $F_n$-equivariant homothety 
between $T,T^\prime$ then 
${\rm CV}(F_n)$ is just the space of equivalence 
classes of points in $cv(F_n)$.
The topology on $cv(F_n)$ descends to 
a natural locally compact topology on ${\rm CV}(F_n)$. The  
group ${\rm Out}(F_n)$ acts
properly discontinuously on ${\rm CV}(F_n)$.  
Write $\partial{\rm CV}(F_n)$ to denote the projectivization
of $\partial cv(F_n)$. Then 
$\overline{{\rm CV}(F_n)}={\rm CV}(F_n)\cup
\partial {\rm CV}(F_n)$ is a compact ${\rm Out}(F_n)$-space.

{\bf Notational convention:} In the sequel we are going to use the
following notations.
\begin{enumerate}
\item An element of $F_n$ is denoted by a small letter, and $[w]$ is the 
conjugacy class of $w\in F_n$.
\item A point in $\overline{cv(F_n)}=cv(F_n)\cup 
\partial cv(F_n)$ is denoted by a capital letter, 
and $[T]\in \overline{{\rm CV}(F_n)}$ 
is the
projectivization of $T\in \overline{cv(F_n)}$.
\item A current is denoted by a Greek letter, and $[\nu]$ is the 
projectivization of $\nu\in {\rm Curr}(F_n)$.
\end{enumerate}

For every $T\in \overline{cv(F_n)}$
and every
$w\in F_n$, the \emph{translation length}
$\Vert w\Vert_T$ for the action of $w$ on $T$ 
is defined to be the \emph{dilation} 
\[\Vert w\Vert_T=
\inf\{d(x,wx)\mid x\in T\}.\]
If $w^\prime\in F_n$ is 
conjugate to $w$ then
$\Vert w^\prime\Vert_T=\Vert w\Vert_T$ and hence
the translation length  $\Vert [w]\Vert_T$ of the conjugacy
class $[w]$ of $w$ is defined.

To every $T\in \partial cv(F_n)$ 
we can associate a topological 
lamination $L(T)$ of zero-length geodesics
\cite{CHL07a} as follows.
For every $\epsilon >0$ define $\Omega_\epsilon(T)$
to be the set 
of all elements $w\in F_n$ with translation length
$\Vert w\Vert_T<\epsilon$. Denote by
$L_\epsilon(T)$ the smallest $F_n$-invariant closed subset
of $\partial^2(F_n)$ which contains all pairs of fixed points of
each element in $\Omega_\epsilon(T)$.
Then 
\[L(T)=\cap_{\epsilon >0}L_\epsilon(T)\]
is a nonempty closed
$F_n$-invariant flip invariant 
subset of $\partial^2 F_n$
which will be called the
\emph{zero lamination} of $T$. Note that two $\mathbb{R}$-trees
$T,T^\prime \in \partial cv(F_n)$ with the same
projectivization have the same
zero lamination. Thus the zero lamination is defined for points in
$\partial {\rm CV}(F_n)$.

The following result is due to Kapovich and Lustig \cite{KL07b,
KL07c}.

\begin{proposition}\label{kaplust}
\begin{enumerate}
\item There is a unique continuous ${\rm Out}(F_n)$-invariant
\emph{length pairing}
\[\langle,\rangle:\overline{cv(F_n)}\times {\rm Curr}(F_n)
\to [0,\infty)\] which satisfies 
$\langle T,\eta\rangle=\Vert [w]\Vert _T$ for every current $\eta$
induced by an indivisible 
conjugacy class $[w]$ in $F_n$ and for every $T\in \overline{cv(F_n)}$.
\item If $T\in \partial cv(F_n)$ then 
$\langle T,\nu \rangle=0$ if and only
if $\nu$ is supported in the zero lamination of $T$.
\end{enumerate}
\end{proposition}

\begin{remark} Kapovich and Lustig call the length
pairing as defined above an intersection form.
\end{remark}

Define ${\cal U\cal M\cal L}^\prime\subset {\cal P\cal M\cal L}(F_n)$ to 
be the set of all projective measured laminations 
$[\nu]$ with the property that $\langle [T|,[\nu]\rangle =0$ for
precisely one projective tree $[T]\in \partial{\rm CV}(F_n)$.
Note that this makes sense without referring to 
specific representatives of the projective classes.
The projective tree $[T]$ is called \emph{dual} to $[\nu]$.
We denote by ${\cal U\cal T}^\prime\subset \partial{\rm CV}(F_n)$ 
the set of all projective trees which are dual to points in 
${\cal U\cal M\cal L}^\prime$.

By invariance of the length pairing, the sets 
${\cal U\cal M\cal L}^\prime$
and ${\cal U\cal T}^\prime$ are invariant under the action of 
${\rm Out}(F_n)$. The assignment $\omega^\prime$ which associates to
$[\nu]\in {\cal U\cal M\cal L}^\prime$ the projective
tree $\omega^\prime([\nu])\in {\cal U\cal T}^\prime$ which is dual
to $[\nu]$ is
${\rm Out}(F_n)$-equivariant.
However, this map is not injective as the 
following example shows. 
This example was provided by the referee of
an earlier version of this paper.

\begin{example}\label{example1}
Let $S$ be a compact surface 
of genus $g\geq 4$ with connected boundary $\partial S$.
The fundamental group of $S$ is the group $F_{2g}$.
Let $\alpha$ be a simple closed separating curve
on $S$ which decomposes $S$ into a surface $S_0$
of genus $2$ and a surface $S_1$ of genus $g-2$ with 
boundary $\partial S\cup\alpha$. Let $\mu$ be a uniquely
ergodic measured geodesic lamination on the surface $S$ 
which is supported in $S_1$ and
is maximal with this property.
This means that after replacing the two 
boundary circles of $S_1$ by punctures, the support of 
$\mu$ decomposes $S_1$ into trigons and two once punctured
monogons. Let $T$ be the tree which is dual
to $\mu$. Since $\mu$ is uniquely ergodic, 
the tree $T$ supports a unique transverse measure up to scaling.

Now let $\nu,\nu^\prime$ be two distinct currents whose
support is the full group $\pi_1(S_0)$. Then the currents
$\mu+\nu,\mu+\nu^\prime$ are both dual to $T$ and to no
other tree. Namely, any tree dual to $\mu+\nu$ contains
in its zero lamination 
the support of $\mu$ as well as the
support of $\nu$. Since the support
of $\nu$ equals $\pi_1(S_0)$ and $\mu$ is uniquely ergodic,
such a tree is topologically equivalent to $T$. But
the tree $T$ supports a unique
transverse measure up to scaling and hence
a tree dual to $\mu+\nu$ is projectively equivalent to $T$.
This reasoning also applies to $\mu+\nu^\prime$ and shows
that a tree which is dual to $\mu+\nu,\mu+\nu^\prime$ is
projectively equivalent to $T$.

To complete the example we have to show that
the currents $\mu+\nu,\mu+\nu^\prime$ are contained in 
${\cal M\cal L}(F_n)$.
For this note that by Proposition \ref{properfreefactor}
and the fact that ${\cal M\cal L}(F_n)$ is a \emph{closed}
subset of ${\rm Curr}(F_n)$, a current of the form
$\zeta+\nu$ is contained in ${\cal M\cal L}(F_n)$ for any
weighted dual current $\zeta$ of a primitive conjugacy
class in $\pi_1(S_1)<\pi_1(S)$. Since $\mu$ is a measured
geodesic lamination on $S_1$, it can be approximated 
in the space of currents by weighted duals of primitive
conjugacy classes in $S_1$. This completes the example.
\end{example}

Let 
\begin{align}
{\cal U\cal M\cal L}=& \{[\nu]\in {\cal U\cal M\cal L}^\prime\mid\notag\\
\langle\omega^\prime[\nu],[\zeta]\rangle =0 & \text{ for }
[\zeta]\in {\cal P\cal M\cal L}(F_n) \text{ only if }
[\zeta]=[\nu]\}.\notag\end{align}
By equivariance, the
set ${\cal U\cal M\cal L}$ is invariant under the action of 
${\rm Out}(F_n)$. The restriction
\[\omega=\omega^\prime\vert {\cal U\cal M\cal L}\] 
of the map 
$\omega^\prime$ is a bijection of ${\cal U\cal M\cal L}$ onto
an ${\rm Out}(F_n)$-invariant subset ${\cal U\cal T}$ of 
${\cal U\cal T}^\prime$.

\bigskip
\begin{remark}
\begin{enumerate} \item
The set ${\cal U\cal M\cal L}$ can be viewed as
the analogue of the set of all projective classes of 
measured geodesic laminations for 
a surface of higher genus whose support fills up $S$ and is
uniquely ergodic (compare Example \ref{example1}).
\item The sets ${\cal U\cal M\cal L}$ and ${\cal U\cal T}$ 
are characterized by having unique duals.
In other words, the definition is symmetric, and we
could begin with defining a set of projective
trees whose zero lamination supports a unique
projective measured lamination etc. The discussion is
completely formal and will be omitted.
\end{enumerate}
\end{remark}

An element $\alpha\in {\rm Out}(F_n)$ is called 
\emph{fully irreducible} (or iwip for short) 
if there is no $k>0$ such that $\alpha^k$ preserves
a free factor of $F_n$. 

The following result describes the action of an iwip element on
the boundary $\partial {\rm CV}(F_n)$ of Outer space 
and on the space of projective measured laminations.
Its first part is due to Levitt and Lustig \cite{LL03}, its second
part is Theorem 36 of \cite{Ma95} which is attributed to Bestvina.

\begin{proposition}\label{northsouth}
\begin{enumerate}
\item An iwip automorphism of $F_n$ acts with north-south dynamics
on the boundary $\partial{\rm CV}(F_n)$ of outer space.
\item An iwip automorphism of $F_n$
acts on ${\cal P\cal M\cal L}(F_n)$
with north-south dynamics.
\end{enumerate}
\end{proposition}

Our first goal is to strengthen this duality between
the action of an iwip element of ${\rm Out}(F_n)$ 
on projective trees and 
projective measured laminations by showing that the
set
${\cal U\cal T}\subset \partial {\rm CV}(F_n)$ 
contains all fixed points of 
iwip elements of ${\rm Out}(F_n)$.
For the proof of this and for later use,
we call an iwip automorphism 
$\alpha\in {\rm Out}(F_n)$ \emph{non-geometric}
if $\alpha$ does not admit any periodic conjugacy class in $F_n$.
This is equivalent to stating that no power of
$\alpha$ can be realized
as a homeomorphism of a compact surface with fundamental group $F_n$
(Theorem 4.1 of \cite{BH92}).

\begin{lemma}\label{uniquefix}
Any fixed point of an
iwip-auto\-mor\-phism on the boundary $\partial {\rm CV}(F_n)$
of Outer space is contained in ${\cal U\cal T}$.
\end{lemma}
\begin{proof} Consider first a non-geometric
iwip-automorphism $\alpha$. 
Let $[T]\in \partial{\rm CV}(F_n)$ 
be the repelling fixed point for the action of
$\alpha$ on the boundary
of Outer space. By Proposition 5.6 of \cite{CHL07b}, 
the zero lamination $L([T])$ of $[T]$ is 
\emph{uniquely ergodic}, i.e. it supports a unique projective
transverse measure $[\nu]$. This projective transverse measure
is a projective measured lamination \cite{Ma95}.
In particular, by the second part of Proposition \ref{kaplust} and
the definitions, 
if $[\nu]\in {\cal U\cal M\cal L}^\prime$ then we also have
$[\nu]\in {\cal U\cal M\cal L}$ and $[T]\in {\cal U\cal T}$.

Now let $\alpha\in {\rm Out}(F_n)$ be a geometric iwip element.
By Theorem 4.1 of \cite{BH92}, 
there is a compact connected surface $S$ with connected
boundary and fundamental group $F_n$ such that 
$\alpha$ can be represented by a pseudo-Anosov
homeomorphism $A$ of $S$.
The repelling projective measured geodesic lamination
$[\nu]$ for $A$ determines up to scaling an action of $F_n=\pi_1(S)$ 
on an $\mathbb{R}$-tree $T$. The 
projectivization $[T]$ of this $F_n$-tree 
is just the repelling fixed point for the action of $\alpha$
on $\partial {\rm CV}(F_n)$. Moreover, $[\nu]$ is supported 
in the zero lamination of $[T]$.

The boundary of $S$ defines a conjugacy class $[w]$ 
in $F_n$ which 
is invariant under $\alpha$.
Any geodesic current supported in the zero lamination 
$L([T])$ of $[T]$ can
be written in the form $a\nu+b\zeta$ where
$\nu$ is a representative of the class $[\nu]$,
where $\zeta$ is the current induced by 
the conjugacy class $[w]$ 
and where $a\geq 0,b\geq 0$. 
We claim that 
if $b>0$ then the 
current $a\nu+b\zeta$ is \emph{not} a 
measured lamination. Namely, otherwise for every
$k>0$ the current $\alpha^k(a\nu+b\zeta)=a\lambda^{-k}\nu+b\zeta$
is a measured lamination as well where $\lambda>1$ is the
expansion rate of $\alpha$. Since the space of measured 
laminations is a closed subset of 
${\rm Curr}(F_n)$ with respect to the weak$^*$-topology, 
this implies that $\zeta\in {\cal M\cal L}(F_n)$.  
However, the invariant conjugacy class 
$[w]$ for $\alpha$ is not
contained in a proper free factor of $F_n$ and hence this
violates Proposition \ref{properfreefactor}.
As a consequence, if $[\nu]\in {\cal U\cal M\cal L}^\prime$ then
also $[\nu]\in {\cal U\cal M\cal L}$ and $[T]\in {\cal U\cal T}$.

We are now left with showing that the projective measured lamination
$[\nu]\in {\cal P\cal M\cal L}(F_n)$ 
defined as above
by the repelling fixed point $[T]\in \partial {\rm CV}(F_n)$ 
of an arbitrary iwip element $\alpha\in {\rm Out}(F_n)$ 
is contained in ${\cal U\cal M\cal L}^\prime$. 
For this let  
$\nu\in {\cal M\cal L}(F_n)$ be a measured lamination
which represents the projective class $[\nu]$.
Assume to the contrary that there is a tree 
$T^\prime\in \partial cv(F_n)$ with 
$[T^\prime]\in \partial {\rm CV}(F_n)-\{[T]\}$ and 
$\langle T^\prime,\nu\rangle =0$.
Since $[T^\prime]\not=[T]$ by assumption and since
by the first part of Proposition \ref{northsouth}
$\alpha$ acts with north-south dynamics on
$\partial{\rm CV}(F_n)$, we have $\alpha^k [T^\prime]\to [Q]$ 
$(k\to \infty)$ where
$[Q]$ is the attracting fixed point for the action of $\alpha$
on $\partial {\rm CV}(F_n)$. 

Choose a continuous section 
$\Sigma:\partial{\rm CV}(F_n)\to \partial cv(F_n)$; such a section
was for example constructed by Skora and White
\cite{S89,W91} (or see the definition below).
Then $\Sigma (\alpha^k [T^\prime])\to \Sigma [Q]$,
moreover we have $\langle \Sigma [Q],\nu\rangle >0$. 
Hence by continuity of the length pairing
on $\partial cv(F_n)\times {\rm Curr}(F_n)$ 
\cite{KL07b} and by naturality under scaling, we conclude
that $\langle\Sigma(\alpha^k[T^\prime]),\nu\rangle>0$
and hence $\langle \alpha^k(\Sigma [T^\prime]),\nu\rangle >0$ 
for sufficiently large $k$. 
Then $\langle \Sigma [T^\prime],\alpha^{-k}\nu\rangle=
\langle \alpha^k(\Sigma [T^\prime]),
\nu\rangle >0$
by invariance of the length pairing
under the action of ${\rm Out}(F_n)$ 
which contradicts the fact that
$\alpha\nu=\rho\nu$ for
some $\rho >0$ and $\langle T^\prime,\nu\rangle =0$. 
\end{proof}

\begin{remark} 
\begin{enumerate}
\item The proof of Lemma \ref{uniquefix}
implies the second part of Proposition \ref{northsouth}.
However, we used Proposition 5.6 of \cite{CHL07b}
which in turn uses a weak version of Martin's result. We also used 
the work of Levitt and Lustig \cite{LL03} which appeared after
Martin's thesis. 
\item In general, there may be points
$[T]\not=[T^\prime]\in \partial {\rm CV}(F_n)$ with the
same zero lamination (see \cite{CHL07c} for a detailed
account on this issue). Lemma \ref{uniquefix} implies
that for the fixed point $[T]$ of an iwip element,
there is no projective tree $[T]\not=[T^\prime]\in \partial {\rm CV}(F_n)$ 
whose zero lamination coincides with the zero lamination of $[T]$.
\end{enumerate}
\end{remark}

We equip ${\cal U\cal M\cal L}$ with the topology as a subspace
of ${\cal P\cal M\cal L}(F_n)$, and we equip ${\cal U\cal T}$ with the
topology as a subspace of
$\partial {\rm CV}(F_n)$.

Denote by 
\[cv_0(F_n)\subset cv(F_n)\]
the subspace of $F_n$-trees $T\in cv(F_n)$ whose
quotient graphs $T/F_n$ have volume one.

For the remainder of this section, choose an arbitrary
tree
\begin{equation}\label{t0}
T_0\in cv_0(F_n)
\end{equation}
and define \begin{equation}\label{Lambda}
\Lambda(T_0)=\{\nu\in {\cal M\cal L}(F_n)\mid
\langle T_0,\nu\rangle=1\}.\end{equation}
We call a lamination $\nu\in \Lambda(T_0)$ \emph{normalized}
for $T_0$, or simply normalized if the 
choice of $T_0$ is clear from the
context.
The next lemma is 
immediate from the continuity of the length pairing. Recall that
the spaces ${\cal M\cal L}(F_n)$ and 
${\cal P\cal M\cal L}(F_n)$ are equipped with the weak$^*$-topology.

\begin{lemma}\label{sectiondependscont}
$\Lambda(T_0)$ is a continuous section of the fibration
\[{\cal M\cal L}(F_n)\to {\cal P\cal M\cal L}(F_n).\]
For each fixed $T_0\in cv_0(F_n)$, the function
$a:cv_0(F_n)\times \Lambda(T_0)\to (0,\infty)$ defined
by $a(T,\zeta)\zeta\in \Lambda(T)$ is continuous.
\end{lemma}

For $T_0\in cv_0(F_n)$ there is a dual section
\begin{equation}\label{Sigma}
\Sigma(T_0)=\{T\in cv(F_n)\cup \partial cv(F_n)\mid 
\max\{\langle T,\nu\rangle\mid
\nu\in \Lambda(T_0)\}=1\}\end{equation}
of the fibration $cv(F_n)\cup \partial cv(F_n)\to {\rm CV}(F_n)\cup 
\partial {\rm CV}(F_n)$.

\begin{remark} A point 
$T\in \Sigma(T_0)\cap cv(F_n)$ can be written in the form
$T=bT^\prime$ where $T^\prime\in cv_0(F_n)$ and where $b>0$
has a geometric meaning. Namely, 
if $d_L$ denotes the one-sided Lipschitz metric on $cv_0(F_n)$ then 
$b=e^{-d_L(T_0,T^\prime)}$. We refer to Corollary \ref{calculatedistance}
for details. 
\end{remark}

For later reference we formulate the
analogue of Lemma \ref{sectiondependscont} 
for the section $\Sigma(T_0)$, but the proof is left to the reader.

\begin{lemma}\label{sectiondependscont2}
$\Sigma(T_0)$ is a continuous section of the fibration
\[cv(F_n)\cup \partial cv(F_n)\to {\rm CV}(F_n)\cup 
\partial {\rm CV}(F_n).\]
For each fixed $T_0\in cv_0(F_n)$ the function
$b:cv_0(F_n)\times \Sigma(T_0)\to (0,\infty)$ defined by
$b(S,T)T\in \Sigma(S)$ is continuous.
\end{lemma}

The following comments summarize some properties
of the sections $\Lambda(T_0)$ and $\Sigma(T_0)$ which will 
be frequently used in the sequel.

\begin{remark}
Let $T_0\in cv_0(F_n)$ be arbitrary.
\begin{enumerate}
\item The spaces
$\Lambda(T_0)$ and $\Sigma(T_0)$ are compact.
\item $0\leq \langle T,\nu\rangle \leq 1$ for all 
$T\in \Sigma(T_0)$ and all $\nu\in \Lambda(T_0)$. 
Moreover, $\langle T,\nu\rangle =0$ only if 
$T\in \partial cv(F_n)$. 
\end{enumerate}
\end{remark}

We use these sections to show

\begin{lemma}\label{dual}
The map $\omega:{\cal U\cal M\cal L}\to 
{\cal U\cal T}$ is an ${\rm Out}(F_n)$-equivariant
homeomorphism.
\end{lemma}
\begin{proof} The map $\omega:{\cal U\cal M\cal L}\to
{\cal U\cal T}$ is clearly 
an ${\rm Out}(F_n)$-equivariant bijection
and hence we just have to show that $\omega$ is 
continuous and open. 

Let $T_0\in cv_0(F_n)$ be a simplicial $F_n$-tree
with quotient of volume one and let
$\Lambda=\Lambda(T_0)$ and $\Sigma=\Sigma(T_0)$ as defined
in equalities (\ref{Lambda}),(\ref{Sigma}). 
By Lemma \ref{sectiondependscont},
the set ${\cal U\cal M\cal L}$ is naturally homeomorphic
to a subset $\Lambda_0$ of $\Lambda$, and ${\cal U\cal T}$ is 
homeomorphic to a subset $\Sigma_0$ of $\Sigma$.
The bijection $\omega$ then induces a bijection
$\omega_0:\Lambda_0\to \Sigma_0$.

Both ${\cal P\cal M\cal L}(F_n)$
and ${\rm CV}(F_n)\cup \partial {\rm CV}(F_n)$ are compact and metrizable
topological 
spaces and hence the same holds true for $\Lambda,\Sigma$. Thus
to show continuity of $\omega$, it suffices to show
that if $(\nu_i)\subset  \Lambda_0$ is any sequence
converging to some $\nu\in \Lambda_0$ then
$\omega_0(\nu_i)\to \omega_0(\nu)$. Via passing to 
a subsequence, we may assume that $\omega_0(\nu_i)\to T$
for some $T\in \Sigma$.
Since $\langle \omega_0(\nu_i),\nu_i\rangle=0$ for all $i$, 
by continuity of the length pairing we have
$\langle T,\nu\rangle=0$ and hence $T=\omega_0(\nu)$.

To show that  $\omega_0$ is open, it suffices to show
that $\omega_0^{-1}:\Sigma_0\to \Lambda_0$
is continuous. However, this follows from 
continuity of the length pairing as in 
the above argument.
\end{proof}

\begin{remark} It can be shown that 
${\cal U\cal M\cal L}$ and ${\cal U\cal T}$ are Borel subsets of
${\cal P\cal M\cal L}(F_n)$ and $\partial{\rm CV}(F_n)$. Since
this fact is not needed in the sequel, we omit a proof. 
\end{remark}

Let $M\subset \partial{\rm CV}(F_n)$ be the closure of the 
set ${\cal U\cal T}$. Since ${\cal U\cal T}$
is invariant under ${\rm Out}(F_n)$, the same holds true for 
$M$. As an easy consequence of Lemma \ref{dual} we obtain:

\begin{corollary}\label{minimalonM}
The action of ${\rm Out}(F_n)$ on $M$ is minimal.
\end{corollary}
\begin{proof}
Lemma \ref{uniquefix} shows that the set ${\cal F}\subset \partial 
{\rm CV}(F_n)$ of fixed points of 
iwip-elements of ${\rm Out}(F_n)$ 
is contained in ${\cal U\cal T}$. This set is invariant 
under the action of ${\rm Out}(F_n)$ and hence
its image under the homeomorphism $\omega^{-1}$ from Lemma \ref{dual}
is invariant as
well. Since the action of ${\rm Out}(F_n)$ on ${\cal P\cal M\cal L}(F_n)$ is 
minimal, Lemma \ref{dual} then implies that ${\cal F}$ is dense in 
${\cal U\cal T}$ and hence in $M$.
  
By the first part of Proposition \ref{northsouth}, 
every iwip element acts with north-south dynamics
on $\partial {\rm CV}(F_n)$. Moreover,  there is no global fixed
point for the action of ${\rm Out}(F_n)$ on $M$. 
Thus the closure of every
${\rm Out}(F_n)$-orbit on $M$ 
contains ${\cal F}$. Since ${\cal F}\subset M$ is 
dense, this shows the lemma.
\end{proof}

\begin{remark} The proof of Corollary \ref{minimalonM} also
implies that a closed ${\rm Out}(F_n)$-invariant
minimal subset of $\partial {\rm CV}(F_n)$ 
is unique. Again, the argument is an immediate
consequence of the work of Levitt and Lustig \cite{LL03}.
Earlier Guirardel \cite{G00} described 
a minimal invariant set for the action of ${\rm Out}(F_n)$ 
explicitly and established uniqueness with a different argument.
\end{remark}

For $\epsilon >0$ define the \emph{$\epsilon$-thick
part} of $cv_0(F_n)$ to be 
the set 
\begin{equation}\label{thick}
{\rm Thick}_\epsilon(F_n)\subset cv_0(F_n)\end{equation}
of all simplicial trees $T\in cv_0(F_n)$
with the additional property that the smallest 
translation length on $T$ of any 
element $w\not=e$ of $F_n$ is at least $\epsilon$.
For sufficiently small $\epsilon>0$ the set
${\rm Thick}_\epsilon(F_n)$  
is a closed connected ${\rm Out}(F_n)$-invariant
subset of $cv_0(F_n)$ on which ${\rm Out}(F_n)$ acts
cocompactly.

The following simple observation will be used
several times in the sequel. For its formulation, recall from
(\ref{Sigma}) the definition of the set $\Sigma(T_0)$ for
some tree $T_0\in {\rm Thick}_\epsilon(F_n)$.

\begin{lemma}\label{scalingzero}
Let $(T_i)\subset {\rm Thick}_\epsilon(F_n)$ be a
sequence such that for every compact set 
$K\subset {\rm Thick}_\epsilon(F_n)$, we have
$T_i\in K$ only for finitely many $i$.
For each $i$ let 
$a_i>0$ be such that
$a_iT_i\in \Sigma(T_0)$. Then $a_i\to 0$ $(i\to \infty)$.
\end{lemma}
\begin{proof}
We follow Kapovich and Lustig \cite{KL07b}.

Let $(T_i)\subset {\rm Thick}_\epsilon(F_n)$ be a sequence
such that for every compact set $K$, we have $T_i\in K$ only 
for finitely many $i$.
For each $i$ let $a_i>0$ be such that 
$a_iT_i\in \Sigma(T_0)$.
Since $\Sigma(T_0)$ is
compact, after passing to a subsequence we may assume that
$a_iT_i\to T$ for some
 $T\in \Sigma(T_0)$. Since $(T_i)$ exits every compact set we have
 $T\in \partial cv(F_n)$.

If $a_i\not\to 0$ $(i\to \infty)$ 
then after passing to a subsequence, we
may assume that $a_i\geq a>0$ for all $i$. Since $T\in \partial cv(F_n)$
there are nontrivial elements in $F_n$ acting on $T$ with 
arbitrarily small translation length. This means that
there exists a sequence $(g_j)\subset F_n-\{e\}$ with 
\begin{equation}\label{limit2}
\lim_{j\to \infty}\Vert g_j\Vert_T=0.
\end{equation}

However, $T_i\in {\rm Thick}_\epsilon(F_n)$ and hence
$\Vert g_j\Vert_{T_i} \geq \epsilon$ for all $i,j$.
Hence for all $i,j$
we have \[a_i\Vert g_j\Vert_{T_i}\geq a\epsilon.\]
On the other hand, for all $j$ 
we have $a_i\Vert g_j\Vert_{T_i}\to \Vert g_j\Vert_T$
$(i\to \infty)$  
which contradicts (\ref{limit2}) above.
Thus indeed $\lim_{i\to \infty}a_i=0$.
This shows the lemma.
\end{proof}

Let $[{\rm Thick}_\epsilon(F_n)]\subset {\rm CV}(F_n)$  be the
projectivization of ${\rm Thick}_\epsilon(F_n)$ and 
denote by  $\overline{[{\rm Thick}_\epsilon(F_n)]}$ the closure of 
$[{\rm Thick}_{\epsilon}(F_n)]$ in 
${\rm CV}(F_n)\cup \partial{CV}(F_n)$. Then 
\[\partial [{\rm Thick}_\epsilon(F_n)]=\overline{[{\rm Thick}_\epsilon(F_n)]}
-[{\rm Thick}_\epsilon(F_n)]\] 
is a closed ${\rm Out}(F_n)$-invariant subset of $\partial{\rm CV}(F_n)$.

For a simplicial tree $T\in {\rm Thick}_\epsilon(F_n)$ call
a primitive conjugacy class $[w]$ in $F_n$ \emph{basic} for $T$
if $[w]$ can be represented by a loop in $T/F_n$ of length at most
two. For example, if $[w]$ can be represented by a loop which
travels through each edge of $T/F_n$ at most twice then
$[w]$ is basic for $T$.
The following duality statement is motivated by 
Teichm\"uller theory.

\begin{lemma}\label{basicclosure}
Let $C\subset \overline{[{\rm Thick}_\epsilon(F_n)]}$
be a closed set. Denote by
$V\subset {\cal P\cal M\cal L}(F_n)$ the set of all 
projective measured laminations which are induced by
a basic primitive conjugacy class for a tree
$S\in {\rm Thick}_\epsilon(F_n)$ with $[S]\in C$ 
and let
$\overline{V}$ be the closure of $V$ in 
${\cal P\cal M\cal L}(F_n)$. Then 
\[\overline{V}-V\subset \{[\mu]\mid \exists\,
[T]\in C\cap \partial [{\rm Thick}_\epsilon(F_n)], 
\langle [T],[\mu]\rangle =0\}.\] 
\end{lemma}
\begin{proof}
Let $C\subset \overline{[{\rm Thick}_\epsilon(F_n)]}$ be a closed set.
Let $V\subset {\cal P\cal M\cal L}(F_n)$ be the set of 
all projective measured laminations which are induced by
a basic primitive conjugacy class for some tree
$T\in {\rm Thick}_\epsilon(F_n)$ with $[T]\in C$. 
Let $([\alpha_i])\subset V$ 
be a sequence which converges to some $[\alpha]\in \overline{V}$.
For every $i\geq 0$ 
let $T_i\in{\rm Thick}_\epsilon(F_n)$ 
be such that $[T_i]\in C$ and that 
$[\alpha_i]$ is induced by a 
primitive conjugacy
class $[w_i]$ in $F_n$ which is 
basic for $T_i$.
After passing to a subsequence we may assume that
$[T_i]\to [T_\infty]\in C$.

Consider first the case 
that $[T_\infty]\in [{\rm Thick}_\epsilon(F_n)]\subset {\rm CV}(F_n)$. 
Let $T_\infty\in cv_0(F_n)$ be the representative
tree with quotient of volume one;
then $T_\infty\in {\rm Thick}_{\epsilon}$ and $T_i\to T_\infty$.
In particular, for sufficiently large $i$ there is a
$3/2$-Lipschitz marked homotopy equivalence
$T_i/F_n\to T_\infty/F_n$.
Thus for sufficiently large $i$ the 
class $[w_i]$ can
be represented by a loop in $T_\infty/F_n$ of length 
at most $3$. However, the
number of conjugacy classes which can be represented by
a loop in $T_\infty/F_n$ of length at most $3$ is finite.
This means that there is a conjugacy class $[w]$ in $F_n$ 
so that $[w_i]=[w]$ for infinitely many $i$. Then 
$[\alpha]=[\alpha_i]$ for at least one $i$ and hence
$[\alpha]\in V$.

In the case that $[T_\infty]\in \partial {\rm CV}(F_n)$ 
let $T_\infty\in \Sigma(T_0)\cap \partial cv(F_n)$ 
be a representative of $[T_\infty]$. For $i\geq 0$ 
let $a_i>0$ be such that 
\[a_iT_i\in \Sigma(T_0).\] 
Then $a_iT_i\to T_\infty$ $(i\to \infty)$, and 
since $T_i\in {\rm Thick}_\epsilon(F_n)$ for all $i$, 
Lemma \ref{scalingzero} implies that $a_i\to 0$
$(i\to \infty)$. 

Let $\alpha_i\in 
\Lambda(T_0)$ be the representative of $[\alpha_i]$.
Then $\alpha_i=b_i \alpha_i^\prime$ where
$\alpha_i^\prime$ is induced by a primitive conjugacy class which is
basic for $T_i$ and $b_i>0$.
Now $T_0\in {\rm Thick}_\epsilon(F_n)$ and hence
$\langle T_0,\alpha_i^\prime\rangle\geq \epsilon$ for all $i$.
This shows
that $b_i\leq 1/\epsilon$ for all $i$.

By compactness, we may assume that
$\alpha_i\to \alpha\in \Lambda(T_0)$ where
$\alpha$ is a representative of the class $[\alpha]$.
Since $\langle T_i,\alpha_i^\prime\rangle \leq 2$ 
by assumption, we have $\langle T_i,\alpha_i\rangle \leq 2/\epsilon$
for all $i$. But 
$\langle a_iT_i,\alpha_i\rangle \to \langle T_\infty,\alpha\rangle$ and 
$a_i\to 0$ $(i\to \infty)$ and therefore 
$\langle T_\infty,\alpha\rangle =0$ by continuity. This shows 
that $[\alpha]$ is supported in the zero lamination of $T_\infty$
and completes the proof of the lemma.
\end{proof}

\begin{corollary}\label{bettercontrol}
Let $C\subset \overline{[{\rm Thick}_\epsilon(F_n)]}$
be a closed set and
let $[T]\in {\cal U\cal T}-C$.
Then $\omega^{-1}([T])\in {\cal U\cal M\cal L}$ 
is not contained in the
closure of the set of all projective
measured laminations which are induced by a basic  
primitive conjugacy class for some tree 
$S\in {\rm Thick}_\epsilon(F_n)$ with 
$[S]\in C$.
\end{corollary}
\begin{proof} If $C\subset \overline{[{\rm Thick}_\epsilon(F_n)]}$
is a closed set and if 
$[T]\in {\cal U\cal T}-C$ then $\omega^{-1}([T])$ is not
supported in the zero lamination of any tree $[S]\in C$.
Together with Lemma \ref{basicclosure}, this shows the
corollary.
\end{proof}

\section{Contracting pairs}

In this section we use the length pairing
to single out sets of  
pairs of distinct points in ${\cal P\cal M\cal L}$ 
whose corresponding lengths functions determine
lines in Outer space with strong contraction properties. 
Examples of such pairs include
all pairs of fixed points of iwip elements in ${\rm Out}(F_n)$.
To define these pairs we establish first some properties
of sums of length functions on Outer space.

Fix some small $\epsilon >0$. For the formulation of the
following definition, recall from (\ref{thick}) in 
Section 2 the definition of 
the set ${\rm Thick}_\epsilon(F_n)$.

\begin{definition}\label{proper}
A family ${\cal F}$ of nonnegative 
functions $\rho$ on 
$cv_0(F_n)$ is called \emph{uniformly proper} if for every
$c>0$ there is a compact subset $A(c)$ of ${\rm Thick}_\epsilon(F_n)$ 
such that 
$\rho^{-1}[0,c]\cap {\rm Thick}_\epsilon(F_n)
\subset A(c)$ for every $\rho\in {\cal F}$.
\end{definition}

Call a pair $(\mu,\nu)\in {\cal M\cal L}(F_n)^2$ \emph{positive}
if the function $T\to \langle T,\nu+\mu\rangle$ is positive
on $cv(F_n)\cup \partial cv(F_n)$.
For $T\in cv_0(F_n)$
recall from Section 2 
the definitions  (\ref{Lambda}), (\ref{Sigma}) 
of the compact sets $\Lambda(T)$ and 
$\Sigma(T)$ before and after Lemma \ref{sectiondependscont}.
These sets are the main tools throughout this section.

\begin{lemma}\label{proper1} 
Let $K\subset {\cal M\cal L}(F_n)\times {\cal M\cal L}(F_n)$ be a 
compact set consisting of positive
pairs. Then the family of functions 
$\{\langle \cdot,\mu+\mu^\prime\rangle\mid (\mu,\mu^\prime)\in K\}$
on $cv_0(F_n)$ is
uniformly proper.
\end{lemma}
\begin{proof} 
Let $K\subset {\cal M\cal L}(F_n)\times {\cal M\cal L}(F_n)$
be as in the lemma. Let $T_0\in {\rm Thick}_\epsilon(F_n)$ and let
$\Sigma=\Sigma(T_0)$.
By assumption, we have
$\langle T,\nu+\nu^\prime\rangle>0$ for every $T\in \Sigma$ and every
$(\nu,\nu^\prime)\in K$.
By continuity of the length pairing and compactness of $\Sigma$ 
and $K$ there is then a number $\delta >0$ such that
$\langle T,\mu+\mu^\prime \rangle\geq \delta$ for every $(\mu,\mu^\prime)\in K$
and every $T\in \Sigma$. 

Let $c>0$ and let $A(c)=\{T\in {\rm Thick}_\epsilon(F_n)\mid
\min\{\langle T,\mu +\mu^\prime \rangle\mid (\mu,\mu^\prime)\in K\}\leq c\}$. 
Then $A(c)$ is a closed subset of ${\rm Thick}_\epsilon(F_n)$.
Our goal is to show that $A(c)$ is compact.

For this assume otherwise. Since ${\rm Thick}_\epsilon(F_n)$ is 
locally compact, there is then a
sequence $(T_i)\subset A(c)$ such that for every compact
set $B\subset {\rm Thick}_\epsilon(F_n)$, 
we have $T_i\in B$ only for finitely many $i$, 
and for each
$i$ there is some $(\mu_i,\mu_i^\prime)\in K$
with 
\begin{equation}\label{lowone}
\langle T_i,\mu_i+\mu_i^\prime\rangle\leq c.
\end{equation}
Let $a_i>0$ be such that $a_iT_i\in \Sigma$. 
Since $\Sigma$ is
compact, after passing to a subsequence we may assume that
\begin{equation}\label{limitone}
\lim_{i\to\infty}a_iT_i=T \notag
\end{equation}
in $\Sigma$.
Moreover, since $K$ is compact, after passing to another subsequence 
we may assume that $(\mu_i,\mu_i^\prime)\to 
(\mu,\mu^\prime)\in K$. Then $\langle T,\mu+\mu^\prime\rangle \geq \delta$.

Lemma \ref{scalingzero} shows that 
$\lim_{i\to \infty}a_i=0$.
On the other hand, since $\langle T,\mu+\mu^\prime\rangle\geq \delta$, 
using once
more continuity of the length function we infer from (\ref{lowone}) that
\[0=\lim_{i\to \infty}a_i\langle T_i,\mu_i+\mu_i^\prime\rangle=
\lim_{i\to \infty}\langle a_iT_i,\mu_i+\mu^\prime_i \rangle=
\langle T,\mu+\mu^\prime \rangle\geq \delta.\] 
This is a contradiction and shows 
that
$\{\langle \cdot,\mu+\mu^\prime\rangle \mid (\mu,\mu^\prime)\in K\}$
is indeed uniformly proper.
\end{proof}

A pair of projective measured laminations
$([\mu],[\nu])\in 
{\cal P\cal M\cal L}(F_n)\times {\cal P\cal M\cal L}(F_n)-\Delta$ 
is called \emph{positive} if for any
representatives $\mu,\nu$ of $[\mu],[\nu]$ 
the
function $\langle\cdot,\mu+\nu\rangle$ on 
$cv(F_n)\cup \partial cv(F_n)$ is positive. 
A tree $T\in cv(F_n)\cup \partial cv(F_n)$ is 
called \emph{balanced} for a positive pair 
$(\mu,\nu)\in {\cal M\cal L}(F_n)^2$ if
$\langle T,\mu\rangle=\langle T,\nu\rangle$. Note that
this only depends on the projective class of $T$.
The set 
\[{\rm Bal}(\mu,\nu)\subset cv(F_n)\cup
\partial cv(F_n)\] of all 
balanced trees for $(\mu,\nu)$ is a closed 
subset of $cv(F_n)\cup \partial cv(F_n)$ which is disjoint
from the set of trees on which either $\mu$ or $\nu$ vanishes.
Let moreover 
\[{\rm Min}_\epsilon(\mu+\nu)\subset {\rm Thick}_\epsilon(F_n)\] be the 
set of all points for which 
the restriction of the function 
$T\to \langle T,\mu+\nu\rangle$ to ${\rm Thick}_\epsilon(F_n)$ 
assumes a minimum. 

\begin{remark} It follows from 
continuity of the function $T\to \langle T,\mu+\nu\rangle$,
local compactness of ${\rm Thick}_\epsilon(F_n)$ and
Lemma \ref{proper1} that the set
${\rm Min}_\epsilon(\mu+\nu)$ is non-empty and compact.
\end{remark}

From now on we fix a sufficiently small $\epsilon >0$. 
The definition of a $B$-contracting pair below depends on this
number $\epsilon$. However, as it will become apparent in the 
proof of Proposition \ref{basicproperty2}, a pair which 
is $B$-contracting for a given choice of $\epsilon$ is 
$B^\prime$-contracting for another choice $\epsilon^\prime$. 
Thus even though the quantitative control of a 
$B$-contracting pair depends on the choice of $\epsilon$,
the large-scale properties do not. 

The reason for fixing a number $\epsilon>0$ for the construction
is two-fold. First, the lack of completeness of the one-sided
Lipschitz metric on Outer space leads to a lack of convexity
of length functions. Second, Teichm\"uller theory indicates that no
information is lost. Namely, there are various ways to describe the 
quality of the axis of a pseudo-Anosov element, and one way to 
do this is to project the axis to a fixed thick part of Teichm\"uller space
and measure the contraction properties of this projection. This is 
exactly what is done here.

\begin{definition}\label{contracting}
For $B>1$, a positive 
pair of points \[([\mu],[\nu|)\in 
{\cal P\cal M\cal L}(F_n)\times {\cal P\cal M\cal L}(F_n)-\Delta\] 
is called \emph{$B$-contracting} if for any pair
$\mu,\nu\in {\cal M\cal L}(F_n)$ of representatives
of $[\mu],[\nu]$ there is some "distinguished" 
$T\in {\rm Min}_\epsilon(\mu+\nu)$ with the following properties.
\begin{enumerate}
\item $\langle T,\mu\rangle/\langle T,\nu\rangle\in [B^{-1},B]$.
\item If $\tilde \mu,\tilde \nu\in \Lambda(T)$ 
are representatives of $[\mu],[\nu]$ then 
$\langle S,\tilde \mu+\tilde \nu\rangle\geq 1/B$ for all 
$S\in \Sigma(T).$
\item Let ${\cal B }(T)\subset \Lambda(T)$ be the set of all
normalized 
measured laminations which are up to scaling induced
by a basic primitive 
conjugacy class for a tree 
$U\in {\rm Bal}(\mu,\nu)\cap {\rm Thick}_\epsilon(F_n)$. 
Then $\langle S,\xi\rangle \geq 1/B$
for every $\xi\in {\cal B}(T)$ and every
tree \[S\in \Sigma(T)\cap
\bigl(\bigcup_{s\in (-\infty,-B)\cup 
(B,\infty)}{\rm Bal}(e^s\mu,e^{-s}\nu)\bigr).\]
\end{enumerate}
\end{definition}

\begin{remark}\label{correctransform}
The above definition is symmetric in $[\mu],[\nu]$. 
Moreover, 
if $([\mu],[\nu])$ is a $B$-contracting pair
then $([\mu],[\nu])$ is $C$-contracting for every $C\geq B$,
\end{remark}

\begin{remark}\label{correctransform2} 
By invariance of the length pairing under the
diagonal action of ${\rm Out}(F_n)$, 
$B$-contracting pairs and their defining data transform
correctly under ${\rm Out}(F_n)$. Thus if
$([\mu],[\nu])\in {\cal P\cal M\cal L}(F_n)^2$
is a contracting pair, if $\mu,\nu\in {\cal M\cal L}(F_n)$ 
is a pair of representatives of $[\mu],[\nu]$ and 
if $T\in {\rm Min}_\epsilon(\mu+\nu)$ is a distinguished 
point for $\mu,\nu$ which has properties (1),(2),(3) stated 
in the definition, then for every
$\phi\in {\rm Out}(F_n)$ the pair $(\phi[\mu],\phi[\nu])
\in {\cal P\cal M\cal L}(F_n)^2$ is $B$-contracting, and
the tree $\phi(T)$ is a distinguished point 
for $\phi(\mu),\phi(\nu)$.
\end{remark}

\begin{remark} In Definition \ref{contracting},
the tree $T$ could be replaced by any tree
in ${\rm Min}_\epsilon(\mu+\nu)$ (however at the expense
of changing the constant $B$ by a controlled amount).
This fact can be derived from the discussion in later
sections, but it will not be used. Moreover, singling out
a special tree $T$ turns out to be convenient for the proofs.
\end{remark}

We call a pair $(\mu,\nu)\in {\cal M\cal L}(F_n)^2$
\emph{$B$-contracting} for some $B>0$ if the
pair $([\mu],[\nu])$ of its 
projectivizations is $B$-contracting.

As in the case of lines of minima in 
Teichm\"uller space, we use $B$-contracting 
pairs to constuct coarsely well defined lines in Outer space.
This construction is carried out in detail in Section 4. 
In the remainder of this section we establish some first properties
of $B$-contracting pairs which shows in particular that
the set of $B$-contracting pairs is not empty and in fact
contains many elements which can be identified explicitly.

\begin{proposition}\label{basicproperty1}
For any $B>0$, the set ${\cal A}(B)$ 
of $B$-contracting pairs is an ${\rm Out}(F_n)$-invariant closed
subset of the space of positive pairs in 
\[{\cal P\cal M\cal L}(F_n)\times {\cal P\cal M\cal L}(F_n)-\Delta.\]
\end{proposition}
\begin{proof} By definition, ${\cal A}(B)$ is ${\rm Out}(F_n)$-invariant.

To show that ${\cal A}(B)$ is a closed subset of the set of all positive
pairs in 
${\cal P\cal M\cal L}(F_n)\times {\cal P\cal M\cal L}(F_n)-\Delta$,
let $([\mu_i],[\nu_i])$ be a sequence of 
$B$-contracting pairs converging to a positive pair 
$([\mu],[\nu])\in {\cal P\cal M\cal L}(F_n)
\times {\cal P\cal M\cal L}(F_n)-\Delta$.

Let $\mu,\nu\in {\cal M\cal L}(F_n)$ be preimages of 
$[\mu],[\nu]$ and 
choose any sequence $(\mu_i,\nu_i)\in {\cal M\cal L}(F_n)^2$ 
of pairs of preimages of $[\mu_i],[\nu_i]$ which 
converges to  $(\mu,\nu)$. 
By Lemma \ref{proper1}, the family of functions
\[{\cal F}=\{\mu_i+\nu_i,\mu+\nu\}\] 
is uniformly proper. 
Thus if $T_i\in {\rm Min}_\epsilon(\mu_i+\nu_i)$ is a point
as in the definition of a $B$-contracting pair 
then up to passing to a subsequence, we may assume that the sequence
$(T_i)$ converges to a point $T\in  
{\rm Thick}_\epsilon(F_n)$. 
By continuity of the length pairing, 
we have $T\in {\rm Min}_\epsilon(\mu+\nu)$ and moreover
$\langle T,\mu\rangle/\langle T,\nu\rangle \in [B^{-1},B]$.

Our goal is now to show that $T$ has 
properties (2) and (3) 
in Definition \ref{contracting}. 


To see property (2), note that by 
Lemma \ref{sectiondependscont},  
if $\alpha\in {\cal M\cal L}(F_n)$ and if $a_i>0,a>0$ is such that
$a_i\alpha\in \Lambda(T_i),a\alpha\in \Lambda(T)$ then $a_i\to a$
since $T_i\to T$. 
In particular, if $\tilde \mu_i,\tilde \nu_i\in \Lambda(T_i)$ are
representatives of $[\mu_i],[\nu_i]$ then $\tilde \mu_i\to \tilde \mu\in 
\Lambda(T)$ and $\tilde \nu_i\to \tilde \nu\in \Lambda(T)$
where $\tilde \mu,\tilde \nu$ are representatives
of $[\mu],[\nu]$. 

On the other hand, 
if $S\in cv_0(F_n)\cup \partial cv_0(F_n)$
and if $b_i>0$ is such that $b_iS\in \Sigma(T_i)$ then 
Lemma \ref{sectiondependscont2} shows that $b_i\to b$ where
$bS\in \Sigma(T)$.  Then 
\[1/B\leq \langle b_iS,\tilde \mu_i+\tilde \nu_i\rangle\text{ and }
\langle b_iS,\tilde \mu_i+\tilde \nu_i\rangle
\to \langle bS,\tilde \mu+\tilde \nu\rangle\] 
by continuity of the
length pairing. This shows property (2) in the definition of a 
$B$-contracting pair.

Now let $\xi\in {\cal B}(T)\subset \Lambda(T)$ be as in the 
third part of the definition for the pair $(\mu,\nu)$.
Then there is 
a tree $U\in {\rm Bal}(\mu,\nu)\cap {\rm Thick}_\epsilon(F_n)$ 
such that up to scaling, $\xi$ is induced by a primitive conjugacy class 
which 
can be represented by a loop on $U/F_n$ of length
at most two. 
Let \[S\in \Sigma(T)\cap {\rm Bal}(e^s\mu,e^{-s}\nu)\text{ for some }
s\in (-\infty,-B)\cap (B,\infty).\] 
We have to show that $\langle S,\xi\rangle \geq 1/B$. 

To see that this is the case, 
let $a_i>0$ be such that
$a_i\xi\in \Lambda(T_i)$; then $a_i\to 1$. 
Let $t_i\in \mathbb{R}$ be such that
$U\in {\rm Bal}(e^{t_i}\mu_i,e^{-t_i}\nu_i)$; then $t_i\to 0$
$(i\to \infty)$. 
There is a sequence $b_i\to 1$ 
and for every sufficiently large $i$ there is a number 
$s_i\in (-\infty,B)\cup (B,\infty)$ such that 
$b_iS\in \Sigma(T_i)\cap {\rm Bal}(e^{s_i+t_i}\mu_i,e^{-s_i-t_i}\nu_i)$.

Now $s_i+t_i\to s$ $(i\to \infty)$ and hence
$s_i+t_i\in (-\infty, -B)\cup (B,\infty)$ for sufficiently large $i$.
By the third requirement in the definition of a $B$-contracting pair we have
$\langle b_iS,a_i\xi\rangle \geq 1/B$
for all sufficiently large $i$ and hence 
$\langle S,\xi\rangle\geq 1/B$ by continuity.
This completes the proof of the proposition.
\end{proof}

\begin{remark} We have ${\cal A}(B)\subset {\cal A}(C)$ for 
$B<C$ and hence $\cup_{B>0}{\cal A}(B)$ is a countable
union of closed subsets of the set of positive pairs. 
\end{remark}

Corollary \ref{thintriangle} and the remark thereafter shows that for 
$([\mu],[\nu])\in {\cal P\cal M\cal L}(F_n)$, being a  
$B$-contracting pair for some $B>0$ is a property of 
the individual projective measured laminations
$[\mu],[\nu]$ rather than of the pair. Once again, 
Teichm\"uller theory shows that there are positive pairs
which are not contracting. Such a pair can be constructed
from a minimal filling measured geodesic lamination on 
a compact surface $S$ with connected boundary which is 
not uniquely ergodic. Since we do not need this fact
we do not discuss it in more detail here.

The following proposition is
the key observation in this paper.

\begin{proposition}\label{basicproperty2}
If $([\nu_+],[\nu_-])\in {\cal U\cal M\cal L}^2$ 
is the pair of fixed points of an iwip element of ${\rm Out}(F_n)$ then
$([\nu_+],[\nu_-])$ is $B$-contracting for some $B>0$.
\end{proposition}
\begin{proof}
Let $\phi\in {\rm Out}(F_n)$ be an iwip element with pair of fixed points
$[\nu_+],[\nu_-]\in {\cal U\cal M\cal L}$.  In particular,
$([\nu_+],[\nu_-])$ is a positive pair. Up to exchanging
$\phi$ and $\phi^{-1}$ 
there are numbers $\lambda_+,\lambda_->1$ such that
for any representatives $\nu_+,\nu_-$
of the classes $[\nu_+],[\nu_-]$ we have  
\[\phi\nu_+=\lambda_+\nu_+,\phi\nu_-=\lambda_-^{-1}\nu_-.\] 

Let $s_0>0,a>0$ be such that
$e^{s_0}=a\lambda_+,e^{-s_0}=a\lambda_-^{-1}$.
Then 
\[{\cal F}=\{f_s:T\to \langle T,e^s\nu_++e^{-s}\nu_-\rangle \mid 
s\in [-s_0,s_0]\}\]
is a set of functions on $cv_0(F_n)$ which is compact
with respect to the topology
of uniform convergence on compact sets.
Lemma \ref{proper1} shows that the set 
\[C\subset {\rm Thick}_\epsilon(F_n)\] of all minima
of the restrictions 
of all functions from the 
collection ${\cal F}$  to ${\rm Thick}_\epsilon(F_n)$
is compact. In particular, there is a number $B_1>0$
such that 
\[\langle S,e^s\nu_+\rangle/\langle S,e^{-s}\nu_-\rangle\in [B_1^{-1},B_1]\]
for all $S\in C$ and all $s\in [-1-s_0,s_0+1]$. 
This shows that for $B=B_1$ and for any $s\in [-s_0,s_0]$, 
the first requirement in 
Definition \ref{contracting} is fullfilled for
for $e^s\nu_+,e^{-s}\nu_-$.



To establish property (2), for 
$T\in {\rm Thick}_\epsilon(F_n)$
let $\tilde \nu_+(T),\tilde \nu_-(T)\in \Lambda(T)$
be the representative of $[\nu_+],[\nu_-]$ contained in $\Lambda(T)$.
Then  
\[g_T:S\to \langle S,\tilde \nu_+(T)+\tilde \nu_-(T)\rangle\] 
is a function on $cv_0(F_n)$ which depends
continuously on $T$. In particular, the family
${\cal G}=\{g_T\mid
T\in C\}$ is compact with respect to the 
compact open topology for continuous functions
on $cv_0(F_n)\cup \partial cv_0(F_n)$.

Now $\Sigma(T)$ 
depends continuously
on $T\in cv_0(F_n)$ and therefore 
\[{\cal T}=\cup_{T\in C}\Sigma(T)\] is a
compact subset of $cv(F_n)\cup \partial cv(F_n)$.
The restriction to ${\cal T}$ of
every $g\in {\cal G}$ 
is positive. 
This implies that 
\[b=\inf\{g_T(S)\mid g_T\in {\cal G},S\in {\cal T}\}>0.\]
As a consequence, for $B=1/b$ and  $s\in [-s_0,s_0]$,
property (2) in Definition \ref{contracting} 
is fullfilled for $e^s\nu_+,e^{-s}\nu_-$

To establish property (3), 
let $T\in {\rm Min}_\epsilon(\nu_++\nu_-)$ and let 
$K\subset \Sigma(T)$ be the compact subset of all 
trees which are balanced for $(\nu_+,\nu_-)$
and whose projectivizations are  
contained in $\overline{[{\rm Thick}_\epsilon(F_n)]}$. 
Let $\Theta(K)\subset \Lambda(T)$ 
be the closure of the set of all normalized measured laminations 
which are up to scaling induced by some basic primitive conjugacy class for
any tree which is the projectivization of an element of $K$.
By Corollary \ref{bettercontrol}, 
we have $\tilde\nu_+,\tilde \nu_-\not\in \Theta(K)$.

Let $T_+,T_-\in \Sigma(T)$ be dual to $[\nu_+],
[\nu_-]$. If $\zeta\in {\cal M\cal L}(F_n)$ is any measured
lamination then 
$\langle T_{\pm},\zeta\rangle =0$ only
if the projective class of $\zeta$ equals $[\nu_\pm]$.
By continuity of the length pairing,
the set of functions 
\[{\cal F}=\{T\to \langle T,\zeta\rangle\mid \zeta\in \Theta(K)\}\]
is compact for the compact open topology on 
the space of continuous functions on $\Sigma(T)$.
As a consequence, their values
on $T_+,T_-$ are bounded from
below by a positive number. 

Since $[T_+],[T_-]\in {\cal U\cal T}$,
the sets \[U(p)=\{[S]\in \overline{[{\rm Thick}_\epsilon(F_n)]}
\mid S\in {\rm Bal}(e^t\nu_+,e^{-t}\nu_-)\text{ for some }
t>p\}\] $(p>0)$ form a 
neighborhood basis for $[T_+]$ 
in $\overline{[{\rm Thick}_\epsilon(F_n)]}$.
This implies that there is some $p>0$ 
and a number $c>0$ such that
the functions from the set ${\cal F}$ are 
bounded from below by a positive number on  
$\tilde U(p)=\{S\in \Sigma(T)\mid [S]\in U(p)\}$. In the same way
we can construct a neighborhood $V(p)$ of 
$[T_-]$ in $\overline{[{\rm Thick}_\epsilon(F_n)]}$ 
so that the values of the functions from 
${\cal F}$ on $\tilde V(p)$ are bounded from
below by a positive number.
As a consequence, for $s\in [-s_0,s_0]$,
property (3) in Definition \ref{contracting} holds true for 
$e^2\nu_+,e^{-s}\nu_-$.

If $s\in \mathbb{R}$ is 
arbitrary then there is some $m\in \mathbb{Z}$ and some
$s_1\in [0,s_0)$ such that $s=ms_0+s_1$.
By the choice of $s_0$, the function
\[T\to \langle T,\phi^m(e^{s_1}\nu_+)+\phi^m(e^{-s_1}\nu_-)\rangle\] 
is a multiple of the function 
$T\to \langle T,e^s\nu_++e^{-s}\nu_-\rangle$. 
Since $\phi^m$ acts on $cv(F_n)\times {\cal M\cal L}(F_n)$ 
diagonally as a homeomorphism preserving the length pairing,
the three properties of a contracting pair for
$e^s\nu_+,e^{-s}\nu_-$ 
follow from the corresponding properties for $e^{s_1}\nu_+,
e^{-s_1}\nu_-$. We refer to Remark \ref{correctransform2} which 
explains how the defining objects of a contracting
pair transform under the action of ${\rm Out}(F_n)$.
\end{proof}

\section{Axes of $B$-contracting pairs}

The goal of this section is to relate $B$-contracting pairs 
to the geometry of ${\rm Out}(F_n)$. For this we first 
equip $cv_0(F_n)$ with an ${\rm Out}(F_n)$-invariant
distance as follows.

For trees $T,T^\prime\in cv_0(F_n)$ let 
$d_L(T,T^\prime)$ be the logarithm of the minimal Lipschitz
constant of a marked homotopy equivalence $T/F_n\to T^\prime/F_n$.
Then 
\[d(T,T^\prime)=d_L(T,T^\prime)+d_L(T^\prime,T)\]
is an ${\rm Out}(F_n)$-invariant distance function on
$cv_0(F_n)$ inducing the original topology \cite{FM08}.
This distance $d$ will be called the \emph{symmetrized Lipschitz distance}
in the sequel, and we call the function $d_L$ the
\emph{one-sided Lipschitz metric}.
The group ${\rm Out}(F_n)$ acts properly, isometrically and
cocompactly on ${\rm Thick}_\epsilon(F_n)$ equipped with
the restriction of $d$. Here ${\rm Thick}_\epsilon(F_n)$
is defined as in (\ref{thick}) of Section 2. Unfortunately, 
the symmetrized 
Lipschitz metric $d$ is \emph{not} a geodesic metric.

As in Section 2, 
call a primitive conjugacy class $[w]$ in $F_n$
\emph{basic} for a simplicial tree 
$T\in {\rm Thick}_\epsilon(F_n)$  if $[w]$ can be 
represented by a loop in $T/F_n$ of length at most two.
The following result is due to Tad White (unpublished;
a published account can be found in the paper \cite{FM08}
of Francaviglia and Martino).

\begin{lemma}\label{simpleloops}
For $T,T^\prime\in cv_0(F_n)$, 
\[d_L(T,T^\prime)=\sup 
\{\log \frac{\langle T^\prime,\alpha\rangle}
{\langle T,\alpha\rangle}\mid \alpha\in {\cal M\cal L}(F_n)\}.\]
The supremum is
attained for a measured lamination $\alpha$ 
which is induced by a basic primitive
conjugacy class for $T$.
\end{lemma}
\begin{proof} Let $T,T^\prime\in cv_0(F_n)$. Clearly for every measured 
lamination $\alpha\in {\cal M\cal L}(F_n)$ we have
\[d_L(T,T^\prime)\geq 
\log \frac{\langle T^\prime,\alpha\rangle}
{\langle T,\alpha\rangle}.\]
On the other hand,
Proposition 3.15 of \cite{FM08} states that
$d_L(T,T^\prime)$ is the minimum of the logarithm of the quotients
$\frac{\langle T^\prime,\alpha\rangle}
{\langle T,\alpha\rangle}$
where $\alpha$ passes through
the set of all currents dual to 
a conjugacy class $[w]$ in $F_n$ of the following form. $[w]$ 
can be represented by a
loop $\gamma$ in $T/F_n$ which  
either is 
simple or defines an embedded bouquet of two circles in $T/F_n$ or
defines two disjointly embedded simple closed curves in 
$T/F_n$ joined by a disjoint
embedded arc traveled through twice in opposite direction.
In particular, the length of $\gamma$ is a most two.

It is well known (see p. 197/198 of \cite{M67}) that 
if $\gamma$ is an embedded loop in $T/F_n$ then $\gamma$ 
represents
a primitive conjugacy class in $F_n$. If 
$\gamma=\gamma_1\gamma_2$ where $\gamma_1,\gamma_2$ are 
two embedded loops which intersect in a single point then
$\gamma$ can be obtained from the primitive element 
$\gamma_1$ by a Nielsen move with the primitive
element $\gamma_2$ and once
again, $\gamma$ is primitive. The third case is completely
analogous.
\end{proof}

For the following observation, recall from (\ref{Lambda})
and (\ref{Sigma}) of Section 2 the definitions of the sets
$\Lambda(T)$, $\Sigma(T)$ for a tree $T\in cv_0(F_n)$.
Lemma \ref{simpleloops} implies

\begin{corollary}\label{calculatedistance}
Let $T,S\in cv_0(F_n)$ and let $b>0$ be such that
$bS\in \Sigma(T)$. If $\langle bS,\nu\rangle \geq 1/B$ for some
$B>0$ and some $\nu\in \Lambda(T)$ then 
\[\log \langle S,\nu\rangle \leq d_L(T,S)\leq
\log\langle S,\nu\rangle +\log B.\] 
\end{corollary}
\begin{proof} Let $T,S\in cv_0(F_n)$, $b>0$ be as in the corollary.
By Lemma \ref{simpleloops} and invariance under scaling, we have
\[d_L(T,S)=\sup\{\log \langle S,\alpha\rangle \mid
\alpha\in \Lambda(T)\}.\]
This implies the left hand side of the inequality.

If $b>0$ is such that $bS\in \Sigma(T)$ then 
$b=e^{-d_L(T,S)}$.  Moreover, if $\nu\in \Lambda(T)$ 
and if $\langle bS,\nu\rangle \geq 1/B$ then 
$\langle S,\nu\rangle \geq e^{d_L(T,S)}/B$. Taking the logarithm 
shows the right hand side of the inequality in the corollary.
\end{proof}

From a $B$-contracting pair we now 
construct
a family of "lines" in the $\epsilon$-thick part ${\rm Thick}_\epsilon(F_n)$
of $cv_0(F_n)$.

\begin{definition}\label{axis} An \emph{axis} 
for a $B$-contracting pair $([\mu],[\nu])$ is a map 
\[\gamma:\mathbb{R}\to 
\cup_{s\in \mathbb{R}}{\rm Min}_\epsilon(e^{s/2}\mu+e^{-s/2}\nu)\]
for some representatives $\mu,\nu$ of 
$[\mu],[\nu]$ such that for every $t\in \mathbb{R}$, 
$\gamma(t)$ is a point
in ${\rm Min}_\epsilon(e^{t/2}\mu+e^{-t/2}\nu)$  which has
the properties (1)-(3) as in 
Definition \ref{contracting} for the pair of representatives
$e^{t/2}\mu,e^{-t/2}\nu\in {\cal M\cal L}(F_n)$ 
of the projective measured laminations $[\mu],[\nu]$.
\end{definition}

We do \emph{not} require the map $\gamma$ to be continuous.

Note that an axis as in Definition \ref{axis} depends on 
choices and hence 
is by no means unique. For a fixed choice of representatives
$\mu,\nu$ of the classes $[\mu],[\nu]$, 
the point $\gamma(t)$ of an axis $\gamma$ is chosen
from a fixed compact subset of ${\rm Thick}_\epsilon(F_n)$.

A different choice $\mu_0,\nu_0$ of 
representatives of the projective measured
laminations $[\mu],[\nu]$ yields the same axes, but with
a parametrization which is modified by a 
translation. Namely,
multiplying $\mu,\nu$ with the same positive
scalar does not change ${\rm Min}_\epsilon(e^s\mu+e^{-s}\nu)$ 
for any $s$. On the other hand, if we replace $\nu$ by $a\nu$ for 
some $a>0$ then we have 
\[{\rm Min}_\epsilon(e^s\mu+e^{-s}(a\nu))=
{\rm Min}_\epsilon(e^t\mu+e^{-t}\nu)\]
where $t=s-\frac{1}{2}\log a$.

For a number $D>0$, 
a \emph{$D$-coarse geodesic} in a metric space $(X,d)$ is a 
(possibly non-continuous) map $\gamma:J\to X$ 
defined on a closed connected subset $J$ of $\mathbb{R}$ 
such that
\[d(\gamma(s),\gamma(t))\in [t-s-D,t-s+D]
\text{ for all }s,t\in J,s\leq t.\]
In particular, if $s<t<u$ then 
$d(\gamma(s),\gamma(u))
\geq d(\gamma(s),\gamma(t))+d(\gamma(t),\gamma(u))-3D$.
The main goal of this section is to relate
axes of $B$-contracting pairs to the symmetrized Lipschitz metric.
The following proposition is more generally true for a
map as in Definition \ref{axis} defined by 
a positive pair $([\mu],[\nu])$ which only has
properties (1) and (2) in Definition \ref{contracting}.

\begin{proposition}\label{coarsegeo}
For all $B>0$ there is a number $\kappa_1(B)>0$ such that
an axis of a $B$-contracting pair is a $\kappa_1$-coarse
geodesic for the symmetrized Lipschitz metric.
\end{proposition}

Since axes of $B$-contracting
pairs are entirely contained in ${\rm Thick}_\epsilon(F_n)$
and since ${\rm Out}(F_n)$ acts properly and cocompactly
on ${\rm Thick}_\epsilon(F_n)$, axes of $B$-contracting pairs 
determine a collection of uniform quasi-geodesics for
${\rm Out}(F_n)$ which only depend on 
the length pairing (note that this is true in spite of the
fact that the symmetrized Lipschitz metric is not geodesic). 
This is in analogy 
to Teichm\"uller space where the construction of lines of minima
only uses convexity of length functions \cite{Ke92}.

The proof of Proposition \ref{coarsegeo} is a consequence of 
two simple observations which are used several times in a sequel.

\begin{lemma}\label{functionevaluate}
Let $(\mu,\nu)\in {\cal M\cal L}(F_n)^2$ be a $B$-contracting pair,
let $s\geq 0$ and let $T\in {\rm Min}_\epsilon(\mu,\nu),
S\in {\rm Min}_\epsilon(e^s\mu,e^{-s}\nu)$
be points which have properties (1) and (2) in 
Definition \ref{contracting}. Then 
\[\log (\langle S,\nu\rangle/\langle T,\nu\rangle)
\geq d_L(T,S)-3\log B-\log 2.\]
\end{lemma}
\begin{proof} Via multiplying 
$\mu,\nu$ with a fixed constant (compare the above discussion) we 
may assume without loss of generality that $\nu\in \Lambda(T)$.
By the first property of a $B$-contracting pair, applied to 
both $T$ and $S$, we have
\[\langle T,\mu\rangle \in [B^{-1},B],\,
\langle S,\mu\rangle \leq e^{-2s}B\langle S,\nu\rangle.\]
In particular, if $\tilde \mu\in \Lambda(T)$ is the normalization of 
$\mu$ at $T$
then $\langle S,\tilde \mu\rangle\leq e^{-2s}B^2\langle S,\nu\rangle$
and hence $\langle S,\nu+\tilde\mu\rangle \leq 2B^2\langle S,\nu\rangle.$

Now let $b>0$ be such that $bS\in \Sigma(T)$. 
The second property of a $B$-contracting pair 
and the above estimate shows that 
$1/B\leq 2B^2 \langle bS,\nu\rangle$ and hence 
from Corollary \ref{calculatedistance} 
we infer that
\[d_L(T,S)\leq \log \langle S,\nu\rangle +3\log B+\log 2\]
as claimed.
\end{proof}

\begin{corollary}\label{minisline}
Let $([\mu],[\nu])\in {\cal P\cal M\cal L}(F_n)^2$ be a 
$B$-contracting pair. Let $\mu,\nu$ be representatives of 
$[\mu],[\nu]$, let $s\geq 0$ and let $T\in {\rm Min}_\epsilon(\mu+\nu),
S\in {\rm Min}_\epsilon(e^s\mu+e^{-s}\nu)$
be points which have the properties (1) and (2) in Definition \ref{contracting}.
Then 
\[2s-2\log B\leq d(T,S)\leq 2s+8\log B+2\log 2 .\] 
\end{corollary}
\begin{proof}
Let $T\in {\rm Min}_\epsilon(\mu+\nu),
S\in {\rm Min}_\epsilon(e^s\mu+e^{-s}\nu)$ for some
$s\geq 0$ be points as in Definition \ref{contracting}.
Assume by multiplying $\mu,\nu$ with a fixed constant
that $\nu\in\Lambda(T)$. Lemma \ref{functionevaluate} and 
Corollary \ref{calculatedistance} show that
\begin{equation}\label{geo1}
d_L(T,S)-3\log B-\log 2\leq \log \langle S,\nu\rangle 
\leq d_L(T,S).\end{equation}
By the first property in the definition of a $B$-contracting pair
we have 
\begin{equation}\label{tmu}
\langle T,\mu\rangle \in [B^{-1},B]\text{ and }
\langle S,\nu\rangle/\langle S,\mu\rangle\in 
[e^{2s}B^{-1},e^{2s}B].\end{equation}

Another application of Lemma \ref{functionevaluate} with the
roles of $T,S,\mu,\nu$ exchanged yields
\begin{align}\label{geo2}
d_L(S,T)-3\log B-\log 2 & \leq 
\log( \langle T,\mu\rangle/\langle S,\mu\rangle)\\ 
 &\leq
2s-\log \langle S,\nu\rangle+2\log B.\notag
\end{align} 
Replacing $-\log \langle S,\nu\rangle$ in 
inequality (\ref{geo2}) by the expression
in inequality (\ref{geo1}) shows that
\[d_L(T,S)+d_L(S,T)\leq 2s+8\log B+2\log 2.\]
On the other hand, Lemma \ref{simpleloops} immediately yields that
\[\log\bigl( \langle S,\nu\rangle \langle T,\mu\rangle/
\langle T,\nu\rangle \langle S,\mu\rangle \bigr) \leq d_L(T,S)+d_L(S,T).\]
Together with the estimate (\ref{tmu}) this
implies the lower bound for $d_L(T,S)+d_L(S,T)$ stated in the 
corollary.
\end{proof}

{\it Proof of Proposition \ref{coarsegeo}:} 
Let $\kappa_1=\kappa_1 (B)=8 \log B+2\log 2$. Let
$([\mu],[\nu])\in {\cal P\cal M\cal L}(F_n)$ be a 
$B$-contracting pair, let $\mu_0,\nu_0$ be representatives
of $[\mu],[\nu]$ and let 
$\gamma:\mathbb{R}\to cv_0(F_n)$ be an axis for 
$([\mu],[\nu])$ as in Definition \ref{axis}. 

Let $t<s$ and write $\mu=e^t\mu_0,
\nu=e^{-t}\nu_0$. 
Then 
\[{\rm Min}_\epsilon(e^{s}\mu_0,e^{-s}\nu_0)=
{\rm Min}_\epsilon(e^{s-t}\mu,e^{-(s-t)}\nu)\]
and therefore 
Corollary \ref{minisline} applied to $\mu,\nu$ and 
$e^{s-t}\mu,e^{-(s-t)}\nu$ shows that 
\[2(s-t) -\kappa\leq d(\gamma(2t),\gamma(2s))
\leq 2(s-t)+\kappa\]
as promised.
\qed

The \emph{Hausdorff distance} between two
closed (not necessarily compact) 
sets $A,B\subset {\rm Thick}_\epsilon(F_n)$ is the  
infimum of the  numbers $r\in [0,\infty]$ 
such that $A$ is contained in the
$r$-neighborhood of $B$ and $B$ is contained in the $r$-neighborhood of $A$.
As an immediate corollary, we obtain for 
$\kappa_1(B)=8\log B + 2\log 2$ the following

\begin{corollary}\label{axiswelldefined}
For every $B$-contracting pair $([\mu],[\nu])$, 
the Hausdorff distance between any two axes of 
$([\mu],[\nu])$ is at most $\kappa_1(B)$.
\end{corollary}
\begin{proof} Let $\gamma_1,\gamma_2$ be any two 
axes of a $B$-contracting pair $([\mu],[\nu])$. 
Assume without loss of generality that $\gamma_1(0)
\in {\rm Min}_\epsilon(\mu+\nu)$, $\gamma_2(0)\in 
{\rm Min}_\epsilon(\mu+\nu)$. 
It suffices to show that $d(\gamma_1(s),\gamma_2(s))\leq \kappa_1(B)$
for all $s\in \mathbb{R}$. However, 
the map $\tilde \gamma$ defined by $\tilde \gamma(s)=\gamma_1(s)$ and
$\tilde \gamma(t)=\gamma_2(t)$ for $t\not=s$ is an axis and hence
the distance estimate  is immediate from Corollary \ref{minisline}
and its proof.
\end{proof}

\begin{remark} In fact, since by Proposition \ref{coarsegeo}, 
axes of $B$-contracting pairs with the parametrization as
in the definition of an axis are uniform coarse geodesics, 
any two of them are uniform fellow travellers after perhaps
changing the parametrization by a translation.
\end{remark}

Recall from Lemma \ref{dual} the definition of the
equivariant homeomorphism $\omega:{\cal U\cal M\cal L}\to 
{\cal U\cal T}$. The next proposition explains that the axis of 
a $B$-contracting pair which determines a pair of distinct
points in
$\partial {\rm CV}(F_n)$ indeed connect these points.
 This is in particular relevant in view of Proposition \ref{basicproperty2}
 and the work of Handel and Mosher \cite{HM06}. For its formulation
 we denote as before by $[T]\in {\rm CV}(F_n)$ the 
 projective class of the tree $T\in cv_0(F_n)$.

\begin{proposition}\label{convergence}
Let $\gamma$ be an axis of a $B$-contracting 
pair $([\mu],[\nu])\in {\cal U\cal M\cal L}^2$.
Then $\lim_{t\to \infty}[\gamma(t)]= \omega([\mu])$ in 
$\overline{{\rm CV}(F_n)}$. 
\end{proposition}
\begin{proof} Let $T=\gamma(0)$ and for $s\geq 0$ let
$\beta(s)>0$ be such that $\beta(s)\gamma(s)\in \Sigma(T)$.
Lemma \ref{scalingzero} shows that $\beta(s)\to 0$
$(s\to \infty)$. 

Choose representatives $\mu,\nu\in \Lambda(T)$ of $[\mu],[\nu]$;
then 
$ \langle \beta(s)\gamma(s),\mu+\nu\rangle \leq 2.$
On the other hand, by the definition of an axis and 
Property (1) in Definition \ref{contracting}, we have
\[\langle \gamma(s),e^{s/2}\mu\rangle/
\langle \gamma(s),e^{-s/2}\nu\rangle \in [B^{-1},B].\]
We conclude that
$\lim_{s\to \infty}\langle \beta(s)\gamma(s),\mu\rangle=0$.
As a consequence, the support of $[\mu]$ is contained in the zero lamination
of every tree $[S]\in \partial {\rm CV}(F_n)$ which 
is an accumulation point of a sequence $[\gamma(s_i)]$ where
$s_i\to \infty$. Since $[\mu]\in {\cal U\cal M\cal L}$, a point in $\partial {\rm CV}(F_n)$
whose zero lamination contains the support of $[\mu]$ is unique. The
proposition follows.
\end{proof}

\section{Axes of contracting pairs are contracting}

The goal of this section is to show that 
an axis of a $B$-contracting pair is contracting for the 
Lipschitz distance. We continue to use all assumptions and
notations from the previous sections.

Recall in particular the definition of an axis 
$\gamma:\mathbb{R}\to {\rm Thick}_\epsilon(F_n)$ 
of a $B$-contracting pair $([\mu],[\nu])\in {\cal P\cal M\cal L}^2$.
It is given by the choice of representatives $\mu,\nu\in {\cal M\cal L}$
of the projective classes $[\mu],[\nu]\in {\cal P\cal M\cal L}$
and for each $t\in \mathbb{R}$  by a point in 
${\rm Min}_\epsilon(e^{t/2}\mu +e^{-t/2}\nu)$ with properties 
(1)-(3) from Definition \ref{contracting}. We call this
point the \emph{distinguished point} 
for the pair $(\mu,\nu)$ in ${\rm Min}_\epsilon(e^{t/2}\mu+
e^{-t/2}\nu)$. When no confusion is possible 
we simply talk about the distinguished point or the 
distinguished point in ${\rm Min}_\epsilon(e^{t/2}\mu+
e^{-t/2}\nu)$.

Let $\gamma$ be an axis of a $B$-contracting pair
$([\mu],[\nu])$. 
There is a coarse projection
$\Pi_{\gamma}:cv_0(F_n)\to \gamma(\mathbb{R})$ 
as follows. For $T\in cv_0(F_n)$ 
choose representatives 
$\mu,\nu$ of the classes $[\mu],[\nu]$ such that
$T\in {\rm Bal}(\mu,\nu)$. Note that the measured laminations 
$\mu,\nu$ are unique up to a common rescaling.
Associate to $T$ the distinguished point 
$\Pi_\gamma(T)\in {\rm Min}_\epsilon(\mu+\nu)\cap \gamma(\mathbb{R})$.

Recall that the map $\gamma:\mathbb{R}\to 
{\rm Thick}_\epsilon(F_n)$ may not be continuous,
and in general, 
the coarse projection $\Pi_\gamma$ is discontinuous as well.
However, it associates to a tree $T$ a unique point
$\Pi_\gamma(T)\in \gamma(\mathbb{R})$. It follows from 
Corollary \ref{minisline}, applied to the special case $s=0$, 
that there is a number $\kappa >0$
only depending on $B$ such that for any other choice 
$\gamma^\prime$ of an axis of $([\mu],[\nu])$ and 
any $T\in cv_0(F_n)$ we have
$d(\Pi_\gamma(T),\Pi_{\gamma^\prime}(T))\leq \kappa$
where as before, $d$ denotes the symmetrized Lipschitz
metric on $cv_0(F_n)$.

The following useful fact gives a first idea about the
nature of the map $\Pi$.

\begin{lemma}\label{almostproj}
$d(\Pi_\gamma(\gamma(t)),\gamma(t))\leq 10\log B
+2\log 2$ for all $t\in \mathbb{R}$.
\end{lemma}
\begin{proof} Let $\mu,\nu\in {\cal M\cal L}$ be representatives
of the projective measured laminations which are used to define
the axis $\gamma$.
By property (1) in Definition \ref{contracting}, we have
\[\langle \gamma(t),e^{t/2}\mu\rangle/\langle \gamma(t),e^{-t/2}\nu\rangle
\in [B^{-1},B]\]
for all $t\in \mathbb{R}$. 
Now by definition, if $U=\Pi_\gamma(\gamma(t))$ then 
$U=\gamma(s)$ where $s\in \mathbb{R}$ is such that
\[\langle \gamma(t),e^{s/2}\mu\rangle=
\langle \gamma(t),e^{-s/2}\nu\rangle.\]
Then $\vert t-s\vert \leq \log B$ and hence the lemma follows from
Corollary \ref{minisline}.
\end{proof}

In the statement of the following proposition, we use both the 
one-sided Lipschitz metric $d_L$ as defined in Section 4
(see in particular Lemma \ref{simpleloops} for the 
property which is most important for our purpose) 
and its symmetrization $d$ which is 
a metric, i.e. which satisfies the triangle inequality.
We call the metric $d$ simply the \emph{Lipschitz metric}
in the sequel and use the notation \emph{one-sided
Lipschitz metric} for $d_L$.

\begin{proposition}\label{contraction}
For every $B>0$ there is a number $\kappa_2=
\kappa_2(B)>0$ with the
following property. Let $([\mu],[\nu])$ be a $B$-contracting
pair, let $\gamma$ be an axis for $([\mu],[\nu])$ and 
let $T\in {\rm Thick}_\epsilon(F_n)$.
\begin{enumerate}
\item If $S\in cv_0(F_n)$ is such that 
$d(\Pi_\gamma(T),\Pi_\gamma(S))\geq \kappa_2$ then 
\[d_L(T,S)\geq 
d_L(T,\Pi_\gamma(T))+d_L(\Pi_\gamma(T),\Pi_\gamma(S))+
d_L(\Pi_\gamma(S),S)-\kappa_2.\]
\item If $S\in {\rm Thick}_\epsilon(F_n)$ is such that
$d(\Pi_\gamma(T),\Pi_\gamma(S))\geq \kappa_2$ then 
\[d(T,S)\geq 
d(T,\Pi_\gamma(T))+d(\Pi_\gamma(T),\Pi_\gamma(S))+
d(\Pi_\gamma(S),S)-\kappa_2.\]
\item 
If $S\in \gamma(\mathbb{R})$ is such that
$d(T,S)\leq \inf_td(T,\gamma(t))+1$ then 
$d(S,\Pi_\gamma(T))\leq \kappa_2$.
\end{enumerate}
\end{proposition}
\begin{proof} 
Let $([\mu],[\nu])\in {\cal P\cal M\cal L}(F_n)^2$ 
be a $B$-contracting pair and let $\mu,\nu\in 
{\cal M\cal L}(F_n)$ be representatives
of $[\mu],[\nu]$.
Let $\gamma$ be an axis of $([\mu],[\nu])$ with  
$\gamma(s)\in {\rm Min}_\epsilon(e^{s/2}\mu+e^{-s/2}\nu)$
and write $\Pi=\Pi_\gamma$.

Let $T\in {\rm Thick}_\epsilon(F_n)$ and assume that
$T\in {\rm Bal}(\mu,\nu)$, i.e. that $\Pi(T)=\gamma(0)$.
Let $S\in cv_0(F_n)$ be such that 
\[d(\Pi(T),\Pi(S))\geq 2B+8\log B+2\log 2.\] 
Let $s\in \mathbb{R}$ be such that
$S\in {\rm Bal}(e^{s/2}\mu,e^{-s/2}\nu)$. By the definition of the
projection $\Pi$ we have $\Pi(S)=\gamma(s)$.
By perhaps exchanging
$[\mu]$ and $[\nu]$ we may assume that $s\geq 0$.

Corollary \ref{minisline} shows that 
\[s-2\log B\leq d(\gamma(0),\gamma(s))\leq 
s+8\log B+2\log 2\]
and hence since we assumed that 
$d(\Pi(T),\Pi(S))=d(\gamma(0),\gamma(s))
\geq 2B+8\log B+2\log 2$ we conclude that
$s\geq 2B$.

Let $\alpha$ be a cycle of maximal dilatation
for a marked homotopy equivalence
$T\to \Pi(T)$
with the smallest Lipschitz constant. 
By Lemma \ref{simpleloops},
we may assume that
$\alpha$ is basic for $T$. 

Let $\xi\in \Lambda(\Pi(T))$ 
be induced by $\alpha$ up to scaling
(here as before, $\Lambda(\Pi(T))$ is defined as
in (\ref{Lambda}) in Section 2).
By property (3) in the
definition of a $B$-contracting pair, the hypotheses
of  Corollary \ref{calculatedistance} are satisfied and give
\[\log \langle S,\xi\rangle \geq d_L(\Pi(T),S)-\log B.\]
On the other hand,
$d_L(T,\Pi(T))=-\log \langle T,\xi\rangle$
and therefore
\begin{equation}\label{geo3}
d_L(T,S)\geq \log \bigl(\langle S,\xi\rangle/\langle T,\xi\rangle\bigr) 
\geq d_L(T,\Pi(T))+d_L(\Pi(T),S)-\log B.\end{equation}

Next we claim that there is a universal constant
$C>0$ so that 
\begin{equation}\label{finalsum}
d_L(\Pi(T),S)\geq d_L(\Pi(T),\Pi(S))+d_L(\Pi(S),S)-C.\end{equation}
Together with the estimate (\ref{geo3}),
this shows the first statement in the proposition.

To establish inequality (\ref{finalsum}), 
assume without loss of generality that
$\nu\in \Lambda(\Pi(T))$. Then we obtain 
from the left inequality of Corollary \ref{calculatedistance} that
\begin{equation}\label{sumup1}
d_L(\Pi(T),S)\geq \log \langle S,\nu\rangle,
\end{equation}
on the other hand also 
\begin{equation}\label{sumup2}
\log\langle S,\nu\rangle =
\log(\frac{\langle S,\nu\rangle}{\langle \Pi(S),\nu\rangle})+
\log\langle \Pi(S),\nu\rangle=
\log \langle \Pi(S),\nu\rangle 
+\log \langle S,\tilde \nu\rangle \end{equation}
where $\tilde \nu=\nu/\langle \Pi(S),\nu\rangle\in \Lambda(\Pi(S))$.

Now $s\geq 0$ and hence
Lemma \ref{functionevaluate} shows that
\begin{equation}\label{sumup4}
\log \langle \Pi(S),\nu\rangle \geq d_L(\Pi(T),\Pi(S))-
3\log B-\log 2.
\end{equation}

By the inequalities (\ref{sumup1}),(\ref{sumup2}),(\ref{sumup4})
we are left with showing that
\[\log \langle S,\tilde \nu\rangle \geq d_L(\Pi(S),S)-\tilde C\]
for a universal constant $\tilde C>0$.

To this end let $\tilde \mu\in \Lambda(\Pi(S))$
be the normalization of $\mu$ at $\Pi(S)$.
By property (1) in Definition \ref{contracting}
and the definition of the map $\Pi$, we have
$\langle S,\tilde \nu\rangle \geq 
\langle S, \tilde \mu\rangle/B$.

Property (2) in Definition \ref{contracting}  
and Corollary \ref{calculatedistance}
imply that 
\begin{equation}\label{sumup3}
\log \langle S,\tilde \nu\rangle 
\geq \log \langle S,\tilde \nu+\tilde \mu\rangle /2B
\geq d_L(\Pi(S),S)-2\log B.
\end{equation}
This is what we wanted to show.

Exchanging the role of $T,S$ and adding the 
resulting
inequalities yields the 
second part of the proposition
(with an adjustment of the additive error) provided
that in addition we have $S\in {\rm Thick}_\epsilon(F_n)$.

To show the third part of the proposition, let
$S\in \gamma(\mathbb{R})$ be such that
$d(T,S)\leq \inf_td(T,\gamma(t))+1$.
Lemma \ref{almostproj} shows that
$d(S,\Pi(S))\leq 12 \log B +2\log 2$. 
Inequality (2) in the statement of the proposition
together with the fact that $\Pi(T)\in \gamma(\mathbb{R})$ 
and hence 
$d(T,\Pi(T))\geq d(T,S)-1$
implies that 
$d(\Pi(T),S)\leq C$ where $C>0$ is a universal
constant. The proposition is proven.
\end{proof}

As an immediate consequence, we obtain the following
contraction property (compare \cite{BFu09}).
For pairs of fixed points of iwip
elements, a version of it was established in \cite{AK08}. 

\begin{corollary}\label{contraction2}
For every $B>0$ there is a number $\kappa_3=\kappa_3(B)>0$
with the following property. Let 
$([\mu],[\nu])$ be a $B$-contracting pair, let 
$\gamma$ be an axis for $([\mu],[\nu])$ and let
$T\in {\rm Thick}_\epsilon(F_n)$.
\begin{enumerate}
\item 
Let $r<\inf\{d_L(T,\gamma(t))\mid t\in \mathbb{R}\}$ 
and let $K=\{S\in cv_0(F_n)\mid d_L(T,S)\leq r\}$.
Then the
diameter of $\Pi_\gamma(K)$ 
with respect to the Lipschitz metric 
does not exceed $\kappa_3$.
\item Let $r<\inf\{d(T,\gamma(t))\mid t\in \mathbb{R}\}$ 
and let $K^\prime=\{S\in {\rm Thick}_\epsilon(F_n)\mid
d(T,S)\leq r\}$. Then the
diameter of $\Pi_\gamma(K)$ with respect to the Lipschitz
metric does not exceed $\kappa_3$.
\end{enumerate}
\end{corollary}
\begin{proof} The first part of the corollary
follows from the first part of Proposition \ref{contraction}
together with the fact that 
$d_L(T,\Pi(T))\geq \inf_t d_L(T,\gamma(t))$.
The second part follows in the same way.
\end{proof}

As in Lemma 3.5 of  \cite{BFu09}, we obtain as a corollary

\begin{corollary}\label{geodesicthin}
For all $B>0,D\geq 0$ there is a number $\kappa_4=\kappa_4(B,D)>0$ with the
following property.
Let $\gamma:\mathbb{R}\to 
{\rm Thick}_\epsilon(F_n)$ be an axis of a $B$-contracting pair.
Let $S\in \gamma(\mathbb{R}),T\in {\rm Thick}_\epsilon(F_n)$ and let
$\rho:[0,m]\to {\rm Thick}_\epsilon(F_n)$ be a $D$-coarse geodesic
for the Lipschitz metric connecting
$T$ to $S$. Then $\rho$ passes through the
$\kappa_4$-neighborhood of $\Pi_\gamma(T)$.
\end{corollary}
\begin{proof} Let $\kappa_2>0$ be as in Proposition \ref{contracting}.
Observe first that there is a number
$\chi>\kappa_2$ with the following property.
Let $\gamma:\mathbb{R}\to {\rm Thick}_\epsilon(F_n)$ be an axis of a 
$B$-contracting pair $([\mu],[\nu])$ and let 
$T\in {\rm Thick}_\epsilon(F_n)$. Let $\Pi=\Pi_\gamma$; 
then 
\begin{equation}\label{special}
d(T,\gamma(s))\geq d(T,\Pi(T))+d(\Pi(T),\gamma(s))-\chi
\end{equation}
for every $s\in \gamma(\mathbb{R})$.
Namely, if $d(\Pi_\gamma(0),\gamma(s))\geq 
\kappa_2+12\log B+2\log 2$
then this is immediate from Proposition \ref{contracting} and
Lemma \ref{almostproj}. Otherwise it is
simply the triangle inequality.

Let $\rho:[0,m]\to {\rm Thick}_\epsilon(F_n)$
be a $D$-coarse geodesic connecting
a point $T=\rho(0)\in {\rm Thick}_\epsilon(F_n)$
to a point $\rho(m)\in \gamma(\mathbb{R})$.
Let $k>0$ be 
the supremum
of all numbers $p>0$ such that the closed $p+1$-ball about
$T$ for the metric $d$ 
is disjoint from $\gamma(\mathbb{R})$. We may
assume that $k>0$. 
Since $\rho$ is a $D$-coarse geodesic, we have 
\[k+3\chi\leq d(\rho(0),\rho(k+D+3\chi))
\leq k+2D+3\chi.\]
Moreover, the third part of Proposition \ref{contracting}
shows that $d(\rho(0),\Pi\rho(0))\leq k+\chi$.

If $d(\Pi\rho(0),\Pi\rho(k+D+3\chi))\leq \chi$
then Proposition \ref{contraction} shows that
\begin{align}
d(\rho(0),\rho(m))& \geq
d(\rho(0),\rho(k+D+3\chi))+
d(\rho(k+D+3\chi),\rho(m))-D\notag \\
& \geq k+3\chi+
d(\rho(k+D+3\chi),\Pi\rho(k+D+3\chi))
\notag\\ &+d(\Pi\rho(0)),\rho(m))-2\chi\notag \\
& \geq d(\rho(0),\Pi\rho(0))+
d(\Pi\rho(0),\rho(m))+\chi\notag
\end{align}
which is impossible.

However, if 
$d(\Pi(\rho(0)),\Pi(\rho(k+D+4\kappa_2))\geq \chi$
then the triangle inequality, the second part of 
Proposition \ref{contraction} and 
the estimate (\ref{special})
show that
\begin{align} d(\rho(0),\rho(m))& \geq 
d(\rho(0),\rho(k+D+3\chi))+d(\rho(k+D+3\chi),\rho(m))-D
\notag\\ &\geq 
d(\rho(0),\Pi\rho(0))+
d(\Pi\rho(0),\Pi\rho(k+D+3\chi))\notag\\
 +2d(\Pi\rho(k+D & +3\chi),\rho(k+D+3\chi))
+d(\Pi\rho(k+D+3\chi),\rho(m))-2\chi-D \notag
\end{align}
which is only possible if 
$d(\Pi\rho(k+D+4\kappa_2),\rho(k+D+4\kappa_2)$
is uniformly bounded. Using once more the estimate (\ref{special}),
then $d(\Pi\rho(0),\Pi\rho(k+D+4\kappa_2)$ is uniformly bounded
as well. The corollary follows.
\end{proof}

Corollary \ref{geodesicthin} can be used to construct from distinct
$B$-contracting pairs new pairs which are $C$-contracting for some
$C>B$. To this end, recall that for a $B$-contracting pair
$([\mu],[\nu])$, for representatives 
$\mu,\nu$ of $[\mu],[\nu]$ and for an axis $\gamma$ of
$([\mu],[\nu])$ such that 
$\gamma(0)\in {\rm Min}_\epsilon(\mu+\nu)$ as in the definition of 
an axis, 
the projection
$\Pi_{\gamma}:cv_0(F_n)\to \gamma(\mathbb{R})$
maps $T\in {\rm Bal}(\mu,\nu)$ to $\gamma(0)$.

If $[\mu]\in {\cal U\cal M\cal L}$
then $\langle [T],[\mu]\rangle =0$ for exactly one
projective tree $[T]\in \partial {\rm CV}(F_n)$ 
(the tree $\omega([\mu])$ with the notations from Lemma \ref{dual}).
As a consequence, if $[\mu],[\nu]\in {\cal U\cal M\cal L}$ then 
the projection $\Pi_{\gamma}$ 
extends to a
projection $cv_0(F_n)\cup \partial {\rm CV}(F_n)-\{\omega([\mu]),
\omega([\nu])\}\to \gamma(\mathbb{R})$.

\begin{lemma}\label{Hausdorff} For every $B>0$ there is a number
$C(B)>0$ with the following property. 
Let $([\mu],[\nu])\in {\cal U\cal M\cal L}$
be a $B$-contracting pair with axis 
$\gamma$. Let $[\zeta]\in {\cal U\cal M\cal L}-\{[\mu],[\nu]\}$ and let 
$T=\Pi_{\gamma}(\omega[\zeta])$. Let $\mu,\nu,\zeta\in 
\Lambda(T)$ be representatives of $[\mu],[\nu],[\zeta]$. Then
for $s>B$ the Hausdorff distance between
${\rm Min}_\epsilon(e^s\mu+e^{-s}\nu)$ and 
${\rm Min}_\epsilon(e^s\mu+e^{-s}\zeta)$
is at most $C(B)$.
\end{lemma}
\begin{proof}  Let $\Pi=\Pi_\gamma$ and let $T=\Pi(\omega[\zeta])$  be as
in the lemma. Let $\mu,\nu,\zeta\in \Lambda(T)$ be representatives of 
$[\mu],[\nu],[\zeta]$. Note that $T\in {\rm Bal}(\mu,\nu)$.
The subset of $\partial {\rm CV}(F_n)$ of all projective
trees which are contained in the closure of 
$[{\rm Thick}_\epsilon(F_n)]$ is closed, non-empty and 
${\rm Out}(F_n)$-invariant. Corollary \ref{minimalonM} then shows that 
this set contains ${\cal U\cal T}$.  In particular, the projective
tree $\omega([\zeta])\in {\rm Bal}(\mu,\nu)$ 
is a limit of a sequence $(T_i)\in 
{\rm Thick}_\epsilon(F_n)\cap {\rm Bal}(\mu,\nu)$. 
Lemma \ref{basicclosure} then shows that $[\zeta]$ is a limit
of a sequence of 
projective measured laminations which are induced by 
basic primitive conjugacy classes for the trees $T_i$.  

Let $\zeta\in \Lambda(T)$ be a representative of $[\zeta]$. 
The third property in the definition of a $B$-contracting pair,
applied to $([\mu],[\nu])$, and continuity then show that 
for every tree $S\in \Sigma(T)\cap 
\bigcup_{t\in (B,\infty)}{\rm Bal}(e^t\mu,e^{-t}\nu)$ we have
$\langle S,\zeta\rangle \geq 1/B$.
 This implies
that 
\begin{equation}\label{mutwo}
\log \langle S,\zeta\rangle \geq d_L(T,S)-\log B\text{ for }
S\in cv_0(F_n)\cap 
\bigcup_{t\in (B,\infty)}{\rm Bal}(e^t\mu,e^{-t}\nu).
\end{equation}
Together with Lemma \ref{functionevaluate} we conclude that there is a
number
$\beta_0=\beta_0(B)>0$ such that
\[\vert \log \langle S,\zeta\rangle -\log \langle S,
\nu\rangle \vert \leq \beta_0\]
for all $s>B$ and every $S\in {\rm Min}(e^s\mu+e^{-s}\nu)$.

As a consequence, there is a number $\beta_1=\beta_1(B)>0$ such that
for $S\in {\rm Min}(e^s\mu+e^{-s}\nu)$ we have
\begin{equation}\label{langlebound}
\langle S,e^s\mu+e^{-s}\nu\rangle/
 \langle S,e^s\mu+e^{-s}\zeta\rangle \in [\beta_1^{-1},\beta_1].\end{equation}

On the other hand, using once more Lemma \ref{functionevaluate}
and the
definition of a $B$-contracting pair, there is a number
$\beta_2=\beta_2(B)>0$ 
such that 
\[\langle U,e^s\mu+e^{-s}\zeta\rangle/\langle S,
e^s\mu+e^{-s}\nu\rangle >\beta_1^2\]
whenever $t\geq s$ and 
$U\in {\rm Bal}(e^t\mu,e^{-t}\zeta)\cap 
{\rm Thick}_\epsilon(F_n)$,  $S\in 
{\rm Min}(e^s\mu+e^{-s}\nu)$ and 
$d_L(S,U)\geq \beta_2$.

Since ${\rm Out}(F_n)$ acts on ${\rm Thick}_\epsilon(F_n)$ cocompactly,
there is a number $\beta_3=\beta_3(B)>0$ such that
$d_L(U,Z)\geq \beta_2$ for all $U,Z\in {\rm Thick}_\epsilon(F_n)$ 
with $d(U,Z)\geq \beta_3$.
This implies that $U\not\in {\rm Min}(e^s\mu+e^{-s}\zeta)$ 
if $U\in {\rm Bal}(e^t\mu,e^{-t}\zeta)\cap 
{\rm Thick}_\epsilon(F_n)$ 
is such that $d(U,{\rm Min}(e^s\mu+e^{-s}\nu))\geq \kappa_3$ where
$\kappa_3>0$ only depends on $B$.

By the definition of a $B$-contracting pair the diameter
of ${\rm Min}(e^s\mu+e^{-s}\nu)$ is uniformly bounded 
independent of $s$. The lemma follows.
\end{proof}

As a consequence we obtain 

\begin{corollary}\label{thintriangle}
Let $([\mu_1],[\nu_1]),
([\mu_2],[\nu_2])\in {\cal U\cal M\cal L}^2$ be $B$-contracting pairs.
If $[\mu_1]\not=[\mu_2]$ then $([\mu_1],[\mu_2])$ is 
$C$-contracting for some $C>0$. 
\end{corollary}
\begin{proof}
Since the pair $([\mu_1],[\mu_2])$ is positive by assumption, for any
representatives $\mu_1,\mu_2$ and all $s<t$ the set
$\cup_{u=s}^t{\rm Min}_\epsilon(e^u\mu_1+e^{-u}\mu_2)$ is compact. 

Now let $\gamma_i$ be an axis of $([\mu_i],[\nu_i])$ $(i=1,2)$ and let
$T=\Pi_{\gamma_1}(\omega[\mu_2])$. Choose representatives
$\mu_1,\nu_1,\mu_2\in \Lambda(T)$ of $[\mu_1],[\nu_1],[\mu_2]$. 
By Lemma \ref{Hausdorff}, for $s>B$ the Hausdorff distance
between ${\rm Min}_\epsilon(e^s\mu_1+e^{-s}\nu_1)$ and 
${\rm Min}_\epsilon(e^s\mu_1+e^{-s}\mu_2)$ is at most $C(B)$.
Moreover, if $t>B$ and if $S\in {\rm Bal}(e^t\mu_1,e^{-t}\nu_1)$ then 
$S\in {\rm Bal}(e^s\mu_1,e^{-s}\mu_2)$ for some $s$ so
that $\vert s-t\vert$ is bounded by a constant only depending on $B$.

Apply this reasoning to the line $\gamma_2$ and $[\mu_1]$. 
If $S=\Pi_{\gamma_2}(\omega[\mu_1])$ and if 
$s_0$ is such that $S\in {\rm Min}_\epsilon(e^{s_0}\mu_2+e^{-s_0}\nu_2)$
then we can compare balancing and projections for the half-line
of $\gamma_2$ which begins at $S$ and corresponds to 
$s\to \infty$.  

Since the non-local
property (3) in Definition \ref{contraction} only depends
on the position of balancing points relative to 
the minimal sets of a positive pair, 
the pair $([\mu_1],[\mu_2])$ satisfies
the requirements Definition \ref{contraction} for some $C>0$. 
\end{proof} 

\begin{remark} In general, the number $C>0$ in Corollary \ref{thintriangle}
depends on $[\mu_1],[\mu_2]$ and not only on $B$. 
This can be seen as follows.

Let $S$ be an oriented surface with connected boundary and 
fundamental group $F_n$. The Teichm\"uller space
of $S$ embeds quasi-isometrically into $cv_0(F_n)$. Let 
$\phi$ be a pseudo-Anosov mapping class; then the $\phi$-invariant
Teichm\"uller geodesic is coarsely an axis for $\phi$ viewed as in iwip
element in ${\rm Out}(F_n)$. In particular, its pair of fixed points 
$([\mu],[\nu])$ in ${\cal P\cal M\cal L}(F_n)$ is a 
$B$-contracting pair for some $B>0$.
Let $\alpha$ be a Dehn twist about a
simple closed curve on $S$. Then for every $k>0$, the pair
$(\alpha^k[\mu],\alpha^k[\nu])$ is the pair of fixed points for
$\alpha^{k}\phi\alpha^{-k}$ and hence it is a $B$-contracting pair
for the same $B$. However, as $k\to \infty$, 
the Teichm\"uller geodesic 
connecting $\alpha^k[\nu]$ to $[\mu]$ has longer and longer subsegments
which are contained in the thin part of Teichm\"uller space.
These segments, however, are compact. 
This implies that as $k\to \infty$, the smallest possible contraction constant
$C_k$ for the pair $([\mu],\alpha^k[\nu])$ tends to infinity.
\end{remark}

\section{Iwip subgroups of ${\rm Out}(F_n)$}

In this section  we use the results obtained so far
to shed some light on subgroups of ${\rm Out}(F_n)$ 
which consist of iwip elements.
We begin with an observation of 
Kapovich and Lustig \cite{KL09b}.
For its formulation, we call a pair $(B_1,B_2)$ of disjoint
closed subsets of 
${{\rm Thick}_{\epsilon}(F_n)}\subset cv_0(F_n)$
 \emph{positive} if
the following holds true. Let $K_i\subset {\cal P\cal M\cal L}(F_n)$
be the closure of the set of all 
projective measured laminations which are induced by 
basic primitive conjugacy classes 
for trees in the set $B_i$. 
Then $([\mu_1],[\mu_2])$ is a positive pair for all 
$[\mu_i]\in K_i$ $(i=1,2)$.

The proof of the second part of the following lemma
was communicated to me by Martin Lustig.

\begin{lemma}\label{iwiprecognition}
Let $(B_1,B_2)\subset {\rm Thick}_\epsilon(F_n)^2$ 
be a positive pair
of closed disjoint sets.
Let $\phi\in {\rm Out}(F_n)$ be such that
\[\phi({{\rm Thick}_\epsilon(F_n)-B_2)}\subset B_1\text{ and }
\phi^{-1}({{\rm Thick}_\epsilon(F_n)-B_1)}\subset B_2.\]
Then $\phi$ is iwip, with attracting fixed point in 
the closure of the projectivization $[B_1]$ of $B_1$ and
repelling fixed point in the closure of the projectivization
$[B_2]$ of $B_2$.
\end{lemma}
\begin{proof} Let $\phi\in {\rm Out}(F_n)$ be as in the lemma.
Then $\phi^kT\in B_1$ for every
tree $T\in {\rm Thick}_\epsilon(F_n)-(B_1\cup B_2)$ and 
all $k\geq 1$ and hence 
the subgroup of ${\rm Out}(F_n)$ generated by $\phi$ is infinite.
Moreover, $\phi(B_1)$ is contained in the
interior of $B_1$, and we have
$\phi^{-1}(B_1)\cup B_2={{\rm Thick}_\epsilon(F_n)}$.

For $i=1,2$ let $K_i\subset {\cal P\cal M\cal L}(F_n)$
be the closure of the set 
of all projective measured laminations 
which are induced by
basic primitive conjugacy classes for trees in  $B_i$.
Then $\phi^{-1}(K_1)\cup K_2$ is a closed non-empty
subset
of ${\cal P\cal M\cal L}(F_n)$
containing all 
projective measured laminations which are induced by 
basic primitive conjugacy classes for trees in 
$\phi^{-1}(B_1)\cup B_2={\rm Thick}_\epsilon(F_n)$. Now the closure of the
set of all projective measured laminations which are induced by
basic primitive conjugacy classes for trees in 
${\rm Thick}_\epsilon(F_n)$ is a closed non-empty 
${\rm Out}(F_n)$-invariant subset of ${\cal P\cal M\cal L}(F_n)$
and hence
${\cal P\cal M\cal L}(F_n)=\phi^{-1}(K_1)\cup K_2$ by minimality
\cite{KL07a}. 
This implies in particular that 
$\phi({\cal P\cal M\cal L}(F_n)-K_2)\subset K_1$ 
and similarly $\phi^{-1}({\cal P\cal M\cal L}(F_n)-K_1)\subset K_2$.

As a consequence, if we define 
\[A_1=\cap_i\phi^iK_1,\,A_2=\cap_i\phi^{-i}K_2\]
then $A_i$ is the intersection of  a nested sequence of 
non-empty compact sets
and hence $A_i\not=\emptyset$. Moreover, 
every periodic point for the action of $\phi$ on ${\cal P\cal M\cal L}(F_n)$ 
is contained in $A_1\cup A_2$. Since both
sets $A_i$ are compact and $\phi$-invariant,  
each of the sets $A_1,A_2$ contains at least one periodic point.
If $[\nu_+]\in A_1,[\nu_-]\in A_2$ 
are such periodic points then 
$([\nu_+],[\nu_-])\in {\cal P\cal M\cal L}(F_n)^2$
is a positive pair by definition of a positive pair $(B_1,B_2)$.

By replacing $\phi$ by $\phi^k$ for some $k\geq 1$ we may
assume that $[\nu_+],[\nu_-]$ are fixed points for $\phi$.
Let $\nu_+,\nu_-\in {\cal M\cal L}(F_n)$ be representatives of 
$[\nu_+],[\nu_-]$.  We claim that up to replacing $\phi$ by $\phi^{-1}$ 
there are numbers $\lambda_+>1,\lambda_->1$ such that
$\phi(\nu_+)=\lambda_+ \nu_+$ and $\phi^{-1}(\nu_-)=\lambda_-\nu_-$.

Namely, if up to exchanging $\phi$ and $\phi^{-1}$ we have
$\phi(\nu_+)=\lambda_+\nu_+$, $\phi(\nu_-)=\lambda_-\nu_-$
for some $\lambda_+\leq 1,\lambda_-\leq 1$ then 
$\phi(\nu_++\nu_-)\leq \nu_++\nu_-$ (as functions on 
$cv_0(F_n)$). Thus if $T\in {\rm Thick}_\epsilon(F_n)$ is arbitrary
and if $\nu_+,\nu_-$ are normalized in such a way that
$\nu_+\in \Lambda(T),\nu_-\in \Lambda(T)$ then 
for every $k\geq 0$ we have 
\[\langle \phi^{-k}T,\nu_++\nu_-\rangle =
\langle T,\phi^k(\nu_++\nu_-)\rangle\leq 2.\]
However, by Lemma \ref{proper1}, the function
$\langle \cdot,\nu_++\nu_-\rangle$ on ${\rm Thick}_\epsilon(F_n)$ 
is proper 
which contradicts the fact that ${\rm Out}(F_n)$ acts property discontinuously
on ${\rm Thick}_\epsilon(F_n)$ and that the order of $\phi$ is infinite.

Thus up to replacing $\phi$ by $\phi^{-1}$ we may assume
the following. If $[\nu_+]\in A_1$ (or $[\nu_-]\in A_2$)
is any
fixed point for $\phi$ and if $\nu_+$ (or $\nu_-$) is 
a representative of $[\nu_+]$ (or $[\nu_-]$) then 
$\phi\nu_{\pm}=\lambda_{\pm}\nu_{\pm}$ where $\lambda_+>1,
\lambda_-<1$. Recall that every fixed point for the
action of $\phi$ on ${\cal P\cal M\cal L}$ is contained in
$A_1\cup A_2$.

Now  assume to the contrary that $\phi$ has a power which
is reducible. For simplicity, assume that $\phi$ is reducible
itself. Then up to conjugation, $\phi$ preserves 
a proper free factor $U$ of $F_n$.

We follow the proof of Proposition 6.1 of \cite{KL09b}. Namely,
choose a train track representative for $\phi$. Assume first that there
is more than one stratum. We may assume that the lowest stratum
represents a proper free factor of $F_n$.

If the transition matrix for the restriction of $\phi$ to this stratum is
irreducible then by the argument in the proof of Proposition 6.1 of 
\cite{KL09b}, the map $\phi$ admits two invariant measured
laminations $\nu_-,\nu_+$ carried by the lowest stratum, 
one contracting and one expanding. Moreover, there is an invariant 
projective tree $[T]$ obtained from the top stratum, with all lower
strata collapsed to become elliptic. By construction, the intersection 
number between a representative $T$ of $[T]$ and $\nu_-,\nu_+$ 
vanishes.

In the case that the transition matrix of the bottom stratum is the identity
there is a primitive element which is fixed which violates the assumption 
on $\phi$.

If there is a single stratum then there are again two cases. 
In the first case, the expanding lamination is contained in a free factor
in which case the argument in the proof of Proposition 6.1 of \cite{KL09b}
applies. Otherwise there is an invariant proper free factor which is
elliptic in the tree. However, as before, in this case there is a primitive element
which is fixed. 

Together we obtain a contradiction to the above discussion.
\end{proof}

As in Section 2, let $M$ be the closure of ${\cal U\cal T}$
in $\partial {\rm CV}(F_n)$. By 
Lemma \ref{minimalonM}, the action of ${\rm Out}(F_n)$
on $M$ is minimal. 
Since for sufficiently small $\epsilon >0$ the set 
\[\partial [{\rm Thick}_\epsilon(F_n)]=
\overline{[{\rm Thick}_\epsilon(F_n)]}-
[{\rm Thick}_\epsilon(F_n)]\] is a closed
non-empty subset
of $\partial {\rm CV}(F_n)$ which is 
invariant under the action of ${\rm Out}(F_n)$,
we have $M\subset \overline{[{\rm Thick}_\epsilon(F_n)]}$.
As in \cite{H09} we can use Lemma \ref{iwiprecognition} to show

\begin{corollary}\label{pairdense}
The set of pairs of fixed points of iwip elements 
is dense in $M\times M$.
\end{corollary}
\begin{proof}
Let $V_1,V_2\subset M$ be open disjoint sets. We have to show
that there is an iwip element with attracting fixed point in $V_1$ and
repelling fixed point in $V_2$.

Since ${\cal U\cal T}$ is dense in $M$,  
by making $V_1,V_2$ smaller we may assume that 
if $([\mu],[\nu])$ is any pair of projective measured laminations
so that $[\mu]$ is supported in the zero lamination of a tree
in $V_1$ and $[\nu]$ is supported in the zero lamination of a tree
in $V_2$ then the pair $([\mu],[\nu])$ is positive.
By continuity of the length pairing and by 
Lemma \ref{basicclosure} we may moreover assume that
there are compact disjoint neighborhoods 
$B_1,B_2\subset \overline{[{\rm Thick}_\epsilon(F_n)]}$
of $V_1,V_2$ with the following property. 
The sets $B_1\cap [{\rm Thick}_\epsilon(F_n)]$,
$B_2\cap [{\rm Thick}_\epsilon(F_n)]$ 
are projectivizations of sets $\tilde B_i\subset {\rm Thick}_\epsilon(F_n)$,
and $(\tilde B_1,\tilde B_2)$ is a positive pair of closed disjoint sets
as defined in the beginning of this section.
We also may assume that there is an iwip-element 
$\phi$ whose fixed points $a,b$ are contained in 
$M-B_1-B_2$.

Since $V_1\subset M$  
is open and the action of ${\rm Out}(F_n)$ on $M$ is minimal
and preserves the set of fixed points of iwip elements,
there is an iwip element 
$u\in {\rm Out}(F_n)$ with attracting fixed point in $V_1$. 
The stabilizer in ${\rm Out}(F_n)$ of 
a fixed point of an iwip element is virtually cyclic \cite{BFH97}
and therefore
we may assume that the repelling fixed point of $u$ is distinct
from $a,b$. Now $u$ acts with north-south dynamics on 
$\partial CV(F_n)$ and hence  
up to perhaps replacing $u$ by a nontrivial power we may assume that
$u\{a,b\}\subset V_1$.
Then $v=u\phi u^{-1}$ is an iwip element with both fixed points in $V_1$.
Similarly, there is an iwip element $w$ with both fixed points in $V_2$.

Via replacing $v,w$ by sufficiently high powers we may assume that
$v(M-V_1)\subset V_1,v^{-1}(M-V_1)\subset V_1$ and 
$w(M-V_2)\subset V_2,w^{-1}(M-V_2)\subset V_2$.
Then $wv(M-V_1)\subset V_2,v^{-1}w^{-1}(M-V_2)\subset V_1$.
Moreover, by perhaps replacing $v,w$ by a suitable power
we may assume that the assumptions in 
Lemma \ref{iwiprecognition} are satisfied for $wv$. 
Then $wv$ is an iwip whose pair
of  fixed points is contained in $V_1\times V_2$.
\end{proof}

We also obtain information on subgroups $\Gamma$ of 
${\rm Out}(F_n)$ which contain at least one iwip element.
For this call iwip elements
$\alpha,\beta\in {\rm Out}(F_n)$ \emph{independent} if
the fixed point sets for the action of 
$\alpha,\beta$ on $\partial{\rm CV}(F_n)$ do not coincide.
By Proposition 2.16 of \cite{BFH97}, the stabilizer in ${\rm Out}(F_n)$ of 
a fixed point of an iwip element is virtually cyclic and hence this
means that the fixed point sets of $\alpha,\beta$ are in fact disjoint.

\begin{proposition}\label{uniformgood}
Let  $\Gamma<{\rm Out}(F_n)$ be a subgroup which contains
an iwip element. If 
$\Gamma$ is not virtually cyclic then 
there are two independent iwip elements $\alpha,\beta\in \Gamma$
with the following properties.
\begin{enumerate}
\item The subgroup $G$ of $\Gamma$ generated by
$\alpha,\beta$ is free and consists of iwip elements.
\item There are infinitely many elements $u_i\in G$ $(i>0)$ 
with fixed points $a_i,b_i\in {\cal U\cal T}$ such that for all $i$ the 
${\rm Out}(F_n)$-orbit of $(a_i,b_i)\in {\cal U\cal T}
\times {\cal U\cal T}-\Delta$ is distinct
from the orbit of $(b_j,a_j)$ $(j>0)$ or $(a_j,b_j)$ $(j\not=i)$. 
\end{enumerate} 
\end{proposition}
\begin{proof} 
Let $\Gamma<{\rm Out}(F_n)$ be a subgroup which contains
an iwip element $\alpha$.
Let $[T_+],[T_-]$ be the fixed points
of $\alpha$ in $\partial{\rm CV}(F_n)$. If the set 
$\{[T_+],[T_-]\}$ is invariant under the action of $\Gamma$ then
by Theorem 2.14 of \cite{BFH97}, the group $\Gamma$ is virtually cyclic.

Thus we may assume that there is some $\gamma\in \Gamma$ with
$\gamma [T_+]\in {\cal U\cal T}-\{[T_+],[T_-]\}$. 
Then $\gamma [T_+],\gamma [T_-]$ are the fixed points
of the iwip element $\zeta=\gamma\circ \alpha\circ \gamma^{-1}$.
By Proposition 2.16 of \cite{BFH97}, 
the fixed point sets of $\alpha,\zeta$ in 
$\partial{\rm CV}(F_n)$ 
are disjoint,
and $\alpha,\zeta$ act with north-south dynamics on the
compact space 
$M\subset\partial {\rm CV}(F_n)$. 

The usual ping-pong lemma, applied to the action of $\alpha,\zeta$
on $M$, implies that  
for sufficiently large $k>0,\ell >0$ the subgroup
$G$ of $\Gamma$ generated by $\alpha^k,\zeta^\ell$ is free.
Lemma \ref{iwiprecognition} shows that 
we may assume that this group consists of iwip automorphisms.
In particular, each non-trivial element of $G$ acts with
north-south-dynamics on $M$, with fixed points contained in 
${\cal U\cal T}$.
In the case that $\alpha,\zeta$ are non-geometric this 
is a consequence of 
the main result of \cite{KL09b}.

To show the second part of the proposition
we have to find infinitely many elements in $G$ which are mutually 
not conjugate in ${\rm Out}(F_n)$ and not conjugate to their inverses. 
For this we follow the argument in the proof of Proposition 5.7
of \cite{H09}. Namely,
let $\Delta\subset M\times M$ be the diagonal. 
Let $(a_+,a_-)\in M\times M-\Delta$ be the 
pair of fixed points of $\alpha\in G$. 
The ${\rm Out}(F_n)$-orbit
of $(a_+,a_-)$ 
is a closed subset of $M\times M-\Delta$ (Theorem 5.3 of 
\cite{BFH97}).
Therefore by the ping-pong construction, there is an independent
iwip element $\beta\in G$ which is not conjugate to $\alpha$ in 
${\rm Out}(F_n)$.  

Let $(a_+,a_-),(b_+,b_-)\in M\times M-\Delta$ be 
the pairs of fixed points 
for the action of $\alpha,\beta$ on $M\subset \partial{\rm CV}(F_n)$.
Since $\alpha,\beta$ are not conjugate in 
${\rm Out}(F_n)$, the ${\rm Out}(F_n)$-orbits of
$(a_+,a_-)$ and $(b_+,b_-)$ are distinct.
This implies that there are open neighborhoods
$U_+,U_-$ of $a_+,a_-$ and $V_+,V_-$ of $b_+,b_-$ such that
the ${\rm Out}(F_n)$-orbit of $(a_+,a_-)$ does not intersect
$V_+\times V_-$ and that the ${\rm Out}(F_n)$-orbit
of $(b_+,b_-)$ does not intersect $U_+\times U_-$.
Via replacing $\alpha,\beta$ by suitable powers we may assume
that 
\[\alpha(M-\overline{U}_-)\subset U_+,\,
\alpha^{-1}(M-\overline{U}_+)\subset U_-,\,
\beta(M-\overline{V}_-)\subset V_+\text{ and }
\beta^{-1}(M-\overline{V}_+)\subset V_-.\]

For numbers $n,m,k,\ell>2$
consider the element
\[f=f_{nmk\ell}=\alpha^n\beta^m\alpha^k\beta^{-\ell}\in G.\]
It satisfies $f(\overline{U}_+)\subset U_+,
f^{-1}(\overline{V}_+)\subset V_+$ and hence
the attracting fixed point of $f$ is contained in $U_+$
and its repelling fixed point is contained in $V_+$.

Since $n>2$, the conjugate
$f_1=\beta^{-1}f \beta$ 
satisfies $f_1(\overline{U}_+)\subset U_+$
and $f_1^{-1}(\overline{U}_-)\subset U_-$, i.e. its
attracting fixed point is contained in $U_+$ and its repelling
fixed point is contained in $U_-$. 
Furthermore, since $m>2$, its conjugate
$f_2=\beta^{-1}\alpha^{-n}f \alpha^n\beta$ has its attracting fixed point in 
$V_+$ and its repelling fixed point in $V_-$, and its conjugate
$f_3=\beta^{-1}f\beta$ has its attracting fixed point in $V_-$ and its
repelling fixed point in $V_+$.

As a consequence, $f$ is conjugate to both an element with
fixed points in $U_+\times U_-$ as well as to an
element with fixed points in $V_+\times V_-$. 
This implies that $f$ is not conjugate to either $\alpha$ or $\beta$.
Moreover, since $\alpha$ and $\beta$ can not both be conjugate
to $\beta^{-1}$, by eventually adjusting
the size of $V_+,V_-$ we may assume that
$f$ is not conjugate to $\beta^{-1}$.

We claim that
via perhaps increasing the values of $n,\ell$ we 
can achieve that $f_{nmk\ell}$ is not conjugate to $f^{-1}$.
Namely, as $n\to \infty$, the fixed points of the conjugate
$\alpha^{-n}f_{(2n)mkl}\alpha^{n}=\alpha^n\beta^m\alpha^k
\beta^{-\ell}\alpha^n$ 
of $f_{(2n)mkl}$ converge to the fixed points of 
$\alpha$. Similarly, the fixed points of the conjugate
$\beta^{-\ell}f_{nmk(2l)}^{-1}\beta^{\ell}=
\beta^\ell \alpha^{-k}\beta^{-m}\alpha^{-n}\beta^\ell$ 
of $f_{nmk(2\ell)}^{-1}$ 
converge as $\ell\to \infty$ to the fixed points of $\beta$.
Thus after possibly conjugating with $\alpha,\beta$, if 
$f_{nmk\ell}$ is conjugate in ${\rm Out}(F_n)$ 
to $f_{nmk\ell}^{-1}$ for all $n,\ell$  
then there is a sequence of pairwise distinct
elements $g_i\in {\rm Out}(F_n)$ which map 
a fixed compact subset $K$ of ${\rm Thick}_\epsilon(F_n)$ 
into a fixed compact subset $W$ of 
${\rm Thick}_\epsilon(F_n)$ 
and such that $g_i(a,b)\to (b_+,b_-)$.
Since ${\rm Out}(F_n)$ acts properly discontinuously
on ${\rm Thick}_\epsilon(F_n)$ this is impossible.

Inductively we can construct in this way a 
sequence of elements of 
$G$ with the properties stated in the second part of the
proposition.
\end{proof}

\begin{remark} In analogy to \cite{FM02} we can define a 
\emph{convex cocompact
subgroup} of ${\rm Out}(F_n)$ to be a hyperbolic group 
$\Gamma<{\rm Out}(F_n)$ with the following property.
There is a number $B>0$ and there is a $\Gamma$-equivariant 
embedding $\rho:\partial \Gamma\to {\cal U\cal M\cal L}$ 
such that for all $\xi\not=\zeta\in \partial \Gamma$ the pair
$(\rho(\xi),\rho(\zeta))$ is $B$-contracting. 

Using the results in this note,
it is not hard to see that for any two independent iwip elements 
$\alpha,\beta\in {\rm Out}(F_n)$ and for all sufficiently large $k>0$,
the subgroup of ${\rm Out}(F_n)$ 
generated by $\alpha^k,\beta^k$ is free, consists of iwip elements
and is convex cocompact in this sense. However we do not pursue the
development of 
a theory of convex cocompact subgroups of ${\rm Out}(F_n)$
here and propose this as an open problem.

About four years after this work was carried out,
Bestvina and Feighn \cite{BF11} and Handel and Mosher
\cite{HM11} showed that there are two analogues of a curve graph for
${\rm Out}(F_n)$ which are hyperbolic.  
Thus as in the case of 
mapping class groups, one can define a convex cocompact subgroup of 
${\rm Out}(F_n)$ to be a subgroup so that the orbit map on one of these
graphs is a quasi-isometry. I think 
it is interesting to explore how this is related to the definition suggested
above. We conjecture that our definition is equivalent to 
a quasi-isometric orbit map for the free factor graph considered
in \cite{BF11}.
\end{remark}

\begin{example}\label{example2}
The mapping class group ${\rm Mod}(S)$ of a surface
of genus $g\geq 1$ with one puncture is the subgroup of 
${\rm Out}(F_{2g})$ preserving the conjugacy class of the puncture.
If $\alpha,\beta$ are two independent pseudo-Anosov elements
in ${\rm Mod}(S)$ then there is some $k>0$ such that 
the subgroup $\Gamma$ of 
${\rm Out}(F_{2g})$ generated by $\alpha^k,\beta^k$ 
satisfies the assumptions in Proposition \ref{uniformgood}. 
Note that we can also assume that $\Gamma$ is a Schottky group 
in ${\rm Mod}(S)$ in the sense of \cite{FM02}.
\end{example}

\section{Second bounded cohomology}

This section is devoted to the proof of the theorem from the introduction.
We continue to use the assumptions and notations from Sections 2-6.
We use the construction in Section 2 and Section 6 of \cite{H08}. 
To this end we first formulate a general sufficient condition for 
infinite dimensional second bounded cohomology 
$H_b^2(\Gamma,\ell^p(\Gamma))$ for a discrete group $\Gamma$ 
of isometries of a proper locally path connected metric space.
This condition was established in \cite{H08} for groups acting
properly on ${\rm CAT}(0)$-spaces. However, the ${\rm CAT}(0)$-property
is never used. The statement in the generality needed for the 
main application 
is as follows.

\begin{theorem}\label{sufficient}
Let $X$ be a proper locally path connected metric space and let
$\Gamma$ be a countable group of isometries acting properly
discontinuously on $X$. Assume that the following conditions are
satsfied.
\begin{enumerate}
\item $X$ admits a compactification by adding a boundary $\partial X$.
The isometric action of $\Gamma$ extends to an action on $\partial X$
by homeomorphisms.
\item There is a free subgroup $G$ of $\Gamma$ with 
two generators so that each element 
$e\not=g\in G$ acts on 
$\partial X$ with north-south dynamics. 
The stabilizer of a pair of fixed points in $\partial X$ 
for each element in $G$ is virtually cyclic.
\item There is some $g\in G$ so that the pair 
$(a,b)\in \partial X\times \partial X$ 
of fixed points of $g$ can be connected by 
a $g$-invariant  $B$-contracting coarse
geodesic $\gamma:\mathbb{R}\to X$, i.e. such that
$\gamma(t)\to a$ $(t\to \infty)$ and
$\gamma(t)\to b$ $(t\to -\infty)$. The Hausdorff distance
between any two such coarse geodesics is at most $B$.
\item There is some $h\not=g\in G$ and there is a 
fixed point $b\in \partial X$ for $g$ such that the 
stabilizer in $\Gamma$ of the pair $(b,hb)$ is trivial.
\end{enumerate}
Then $H_b^2(\Gamma,\ell^p(\Gamma))$ is infinite dimensional
for every $p\in (1,\infty)$.
\end{theorem}

Passing from geodesics in \cite{H08} to coarse geodesics
as needed for the application to ${\rm Out}(F_n)$ uses 
Corollary \ref{geodesicthin} which is an immediate
consequence of the
$B$-contraction property in Proposition \ref{contracting}, 
but does not use any property specific
to the situation at hand. Moreover, suitable versions of  
Lemma \ref{Hausdorff} and
Corollary \ref{thintriangle} are also used whose general
version can be established with the arguments given in 
Section 4.

Now let $\Gamma<{\rm Out}(F_n)$ be a subgroup which is not virtually 
abelian and contains an iwip element. 
Our goal is to apply Theorem \ref{sufficient} to the action of 
$\Gamma$ on $X={\rm Thick}_\epsilon(F_n)$. We  
define $\partial X$
to be the complement of the projectivization 
$[X]$ of $X$ in the closure of $[{\rm Thick}_\epsilon(F_n)]$
in ${\rm CV}(F_n)\cup \partial {\rm CV}(F_n)$. 
The group $\Gamma$ acts on $X$ properly
discontinuously, and this
action extends to an action on $\partial X$.
The existence of a subgroup $G<\Gamma$ with 
properties (2) and (3) above was shown in Proposition
\ref{uniformgood} and Proposition \ref{contraction}.
The second part of property 3) follows from 
Lemma \ref{Hausdorff}.

Let $g\in G$ be an element with 
attracting fixed point $[\mu]\in {\cal U\cal M\cal L}$,
repelling fixed point $[\nu]\in {\cal U\cal M\cal L}$ and 
such that the ${\rm Out}(F_n)$-orbit of $([\mu],[\nu])$ is distinct from the
${\rm Out}(F_n)$-orbit of $([\nu],[\mu])$. The existence of such an element is
guaranteed by Proposition \ref{uniformgood} and Lemma \ref{dual}.

\begin{lemma}\label{trivialstab}
There is some $h\in G$ such that the stabilizer of the pair
$(h[\nu],[\nu])$ in ${\rm Out}(F_n)$ is trivial.
\end{lemma}
\begin{proof} Using the notation from Lemma \ref{dual}, 
by Theorem 2.14 of \cite{BFH97} the stabilizer
of the projective tree $\omega([\nu])\in {\cal U\cal T}$ is virtually cyclic.
Since $G$ is not virtually cyclic, there is some $h\in G$ such that
the stabilizer in ${\rm Out}(F_n)$ of the pair $(h\omega([\nu]),\omega([\nu]))$
is trivial. Together with Lemma \ref{dual}, this shows the lemma.
\end{proof}

As a consequence, the action of $\Gamma$ on ${\rm Thick}_\epsilon(F_n)$
has all the properties stated in Theorem \ref{sufficient} and hence the main
theorem from the introduction holds true for $\Gamma$.
  
Nonetheless we discuss a few more details 
of the proof of Theorem \ref{sufficient}
for subgroups $\Gamma$ of ${\rm Out}(F_n)$. Namely,
let ${\cal A}([\nu])$ be the union of the ordered pairs of 
distinct points in $\Gamma[\nu]$ with the $\Gamma$-translates of 
$([\mu],[\nu]),([\nu],[\mu])$. The group $\Gamma$ naturally acts
on ${\cal A}([\nu])$ from the left. Moreover, ${\cal A}([\nu])$ is
contained in ${\cal U\cal M\cal L}^2$ and hence
it can be identified with a subset of the set of pairs of points in the
boundary 
$\partial {\rm CV}(F_n)$ of Outer space.
For $X={\rm Thick}_\epsilon(F_n)$,
the group $\Gamma$ acts on ${\cal A}([\nu])\times X$.
For a suitable choice of a point $T\in X$,
the group $\Gamma$ acts freely on 
${\cal A}([\nu])\times \Gamma T$. In particular,
a $\Gamma$-orbit under this action can be identified with $\Gamma$.

The goal is now to construct a distance $\delta$ on 
$V={\cal A}([\nu])\times \Gamma T$
which is invariant under the action of $\Gamma$ and under the action of
the flip $\iota$ exchanging the two components of the pair
in ${\cal A}([\nu])$. 
The construction is done in such a way that 
functions which are anti-invariant under $\iota$ and H\"older continuous
with respect to this  distance give rise to bounded 
$\ell^p(\Gamma)$-valued two-cocycles 
for $\Gamma$, i.e. $\Gamma$-invariant bounded maps
$\Gamma\times \Gamma\times \Gamma\to \ell^p(\Gamma)$ 
which satisfy the cocycle equation.
For specific choices of these maps, these
cocycles are then shown to define nontrivial bounded cohomology classes.

If there is a number $B>0$ such that 
for all $g,u\in \Gamma$ with $g[\nu]\not=u[\nu]$ the
pair $(g[\nu],u[\nu])$ is $B$-contracting
then coarse geodesic triangles with
endpoints in $G[\nu]\subset {\cal U\cal M\cal L}$ 
are uniformly thin and the construction of a
distance on ${\cal A}([\nu])\times X$ 
parallels the construction of the family of
Gromov distances on the boundary of a hyperbolic  
geodesic metric space. However, in general
the pair $(g[\nu],u[\nu])$ is $C$-contracting for a 
number $C>0$ depending on the pair, and 
only segments of an axis near is endpoints 
is $B^\prime$-contracting for a number $B^\prime$
which is controlled by $B$. The 
main idea in \cite{H08} is to truncate the 
coarse geodesics to achieve our goal.

{\bf Acknowledgements:} The major part of the 
work in this paper
was carried out in fall 2007 while I visited the MSRI
in Berkeley as a participant of the special semester in
Teichm\"uller theory and Kleinian groups. 
The paper arose from inspiring and helpful discussions
with Mladen Bestvina, Ilya Kapovich, Gilbert Levitt and 
Lee Mosher. Thanks to Martin Lustig for
providing the argument in the second part of the proof of
Lemma \ref{iwiprecognition}.
I am in particular grateful to the anonymous referee for 
carefully reading the first version of this manuscript and for 
many helpful comments.

\noindent
MATHEMATISCHES INSTITUT DER UNIVERSIT\"AT BONN\\
ENDENICHER ALLEE 60,\\
53115 BONN, GERMANY\\
e-mail: ursula@math.uni-bonn.de

\end{document}